%% file: arxiv_final.tex
\documentclass[letterpaper, journal]{IEEEtran}
\input{header.tex}
\title{Geometry of Injection Regions\\ of Power Networks}
\begin{document}
\author{\IEEEauthorblockN{Baosen Zhang and David Tse}
\IEEEauthorblockA{EECS Dept. University of California, Berkeley, CA, USA \\
Email: \{zhangbao, dtse\}@eecs.berkeley.edu}}
\maketitle
\begin{abstract}
We investigate the constraints on power flow in networks and its implications to the optimal power flow problem. The constraints are described by the injection region of a network; this is the set of all
vectors of power injections, one at each bus, that can be achieved
while satisfying the network and operation constraints. If there are
no operation constraints, we show the injection region of a network is the set of all injections
satisfying the conservation of energy. If the network has a tree
topology, e.g., a distribution network, we show that under voltage magnitude,
line loss constraints, line flow constraints and certain bus real and reactive power
constraints, the injection region and its convex hull have the same Pareto-front. The Pareto-front is of interest since these
are the the optimal solutions to the minimization of increasing
functions over the injection region. For non-tree networks, we
obtain a weaker result by characterize the convex hull of the
voltage constraint injection region for lossless cycles and certain combinations of cycles and trees.
%%In this case, the injection region with all three types of constraints is convex since it is a intersection of convex sets.
\end{abstract}

\section{Introduction}
Optimal power flow is a classic problem in power engineering. It is
usually given as a static subproblem of the security constraint unit
commitment problem, in the sense that all the network dynamics such
as transients and generator behaviors are abstracted away
\cite{Stott87}. The objective of the optimal power flow problem is
to minimize the cost of power generation in a electrical network
while satisfying a set of operation constraints. The cost functions
are generally taken to be convex and increasing. This problem has
received considerable attention since the late 1960's
\cite{Dommel68}, and many different algorithms have been developed
for it. For a comprehensive review the reader can consult
\cite{Amarnath11} and the references within. Despite all the
efforts, the optimal power flow problem still remains difficult
\cite{Overbye04}.

The optimal power flow problem is difficult for two reasons.
Firstly, the optimization problem is nonlinear since the power
injected at each of the buses in the network depends quadratically
on the voltages at the buses. Secondly, there is typically a large
number of different types of constraints. For example, each bus
might have voltage magnitude together with real and reactive power limits, and each transmission line might have thermal constraints and line flow constraints. Due
to these two reasons, the optimal power flow problem is a non-convex
optimization problem with many constraints, and is therefore
challenging to solve. The traditional approach is to tackle the
problem using various heuristics and approximations. One widely used
method is to use the so called DC flow approximation
%\footnote{The DC
%flow approximation name comes from the fact that historically,
%engineers used analogue devices that approximated power flow with
%current flow. A perhaps more accurate name would simply be small
%angle approximation, but the term DC flow approximation is standard
%and we follow it here.} 
where all the lines are assumed to be
lossless, all voltage magnitude are assumed to be fixed, and all
angle differences are assumed to be small \cite{Stott09}. To
contrast with the DC flow approximation, the original optimal power
flow problem is sometimes called the AC problem. As pointed out in \cite{Stott09}, the DC approximation performs badly if it is not used in conjunction with a full AC solution (so called hot start DC) or if the resistance to inductance (R/X) ratio of the lines are high. To solve the full AC problem, many global optimization heuristics like genetic algorithms are used, and their effectiveness is generally gauged by simulations. But these algorithms do not offer any guarantees about performance and do not offer intuition into the structure of the optimization problem.

A new approach to the traditional optimization methods was taken by
the authors in \cite{Lavaei11}. They made the surprising empirical
observation that in many of the IEEE benchmark networks the optimal power flow problem has the same
optimal value as its convex dual. The main theoretical result is
that for a {\em purely resistive} network and quadratic cost
functions with  positive coefficients, this convex relaxation is
tight. In addition, the result still holds if the purely resistive
network is perturbed by adding a small reactive part. From this and
their observations about the IEEE benchmarks, \cite{Lavaei11}
conjectured that the convex relaxation of the optimal power flow problem is always tight for general networks. Unfortunately this conjecture is not true
since there exist many counter examples \cite{Makarov08,Lesieutre11}. A natural question arises: if the relaxation is not tight in general, is it tight for some specific class of networks? The results \cite{Lavaei11}
showed that for 'almost' purely resistive networks the problem is
convex, but these networks are somewhat unrealistic since practical
power networks are mostly reactive instead of resistive. An impetus
for this paper is to look for some more realistic classes of network
for which the optimal power flow problem is convexified.

One increasingly important class of networks is the distribution network. The electricity network is made up of two layers: the transmission network and the distribution network. The transmission network consists of high voltage lines that connect big generators to cities and towns. The distribution network usually consists of a feeder connected to the transmission network, and low voltage lines that connect to the end consumers. In addition to the line voltages, the two types of networks have different topologies. The transmission network is sparse, but irregular, whereas the distribution network is configured to be a {\em tree} at any one time of operation. Traditionally, the optimal power flow problem is only solved in the transmission network, since the demands in the distribution network are fixed and there is very little generation, so there is nothing to optimize. But this is expected to change significantly under the new 'smart grid' operating paradigm, where demand response and distributed renewable energy will play a predominant role. In the widely discussed demand response mechanism, the demands in the distribution network are decision variables (subjected to some constraints) \cite{Albadi07,DOE06}. Also, due to increased renewable penetration at the demand level (e.g. rooftop solar) and  increased distributed generation, solving the optimal power flow in the distribution network is a legitimate problem and could contribute to various pricing and control operations. For example, we show that the voltage control problem \cite{Oapos04} can be formulated into such a framework. Since the resistance to inductance ($R/X$) ratio is much higher in the distribution network compared to the transmission network, DC approximations would perform poorly. Therefore, the full AC optimal power flow on the distribution network needs to be solved and we show the tree topology of the distribution network simplifies the problem significantly and allows the full AC problem to be efficiently solved in many situations.  

To find out if the optimal power flow problem is convex for a network, we focus on the {\it feasible injection region} of a power network since it allows one to think
about power flow in a more abstract way and is quite useful in
understanding the structure of the problem. The feasible injection
region is simply the feasibility region of the optimal power flow
problem, i.e. the set of all vectors of feasible {\it real power}
injections (both generations and withdraws) at the various buses
that satisfy the given network and operation constraints (including reactive power constraints). For notational convenience, we drop the word feasible and refer to the region as the injection region. Since the optimization problem is solved over the injection region, it is useful to understand the geometry of the region. We model the reactive powers in the network as constraints at the buses. Therefore the injection region is in terms of the real powers, while possibly satisfying some bus reactive power constraints. 

%Analogously, one could define the {\em reactive} injection region of a network. However, here we make an distinction of the real and the reactive power injection regions. This is because in most practical settings, the objective function of the optimization problem is in terms of {\em real} powers only. For example, the cost curve for an generator only includes the real power output; also, the consumers are only charged based on the amount of real power they consume (watt-hours). Therefore the optimization searches over the real injection region to find the optimal solutions. For the reactive powers, we model them as constraints in the problem. That is, since each bus can only supply reactive powers in a certain range, we model it as an constraint at each bus of the networks.  

Unfortunately, the injection region is not convex in general \cite{Lesieutre05}.  Even
though the region is not convex, it still has some desirable
properties for optimization. A subset of the injection region of
particular interest is the {\em Pareto-front}.\footnote{A point in a
set is called Pareto-optimal if any coordinate cannot be decreased
further without increasing at least one other coordinate; the Pareto
front of a set is simply the set of all Pareto-optimal points.} When
minimizing an increasing function over a set, the optimal solutions
are on the Pareto-front.  Therefore, even though the injection
region is not convex, if its Pareto-front is the same as that of its convex hull, the
optimization problem is still easy.

The use of injection region is also useful since it decouples the
optimization problem from the physics of power flow, thus allowing
us to have a higher level view that is often beneficial for other
problems in optimization, control and pricing in power systems. For
example, \cite{Hogan92} showed there is revenue adequacy in the
financial transmission rights markets if the injection region has a
convex Pareto-front. A similar observation is made by
\cite{Lavaei11b} in the context of economic dispatch. This result
then can be used if the DC flow assumption is made or if the network
is such that the AC injection region where the above condition is
true. This is similarly the case for many of the recently proposed
demand response algorithms.

As a starting point, we look at the injection region of a network
with no constraints. In this case, we show the injection region is
simply the upper half space that satisfies the law of conservation
of energy.
%\footnote{This settles an open question in
%\cite{Jarjis81,Hogan92}.} 
Therefore, the difficult and interesting
part is to quantify how the injection region changes once the
operation constraints are added.

There are typically four types of operation constraints in a power
network: voltage magnitude, thermal loss in transmission lines, line
flow limits in a transmission line and bus real and reactive power limits. If the
network is a {\em tree}, we show that under
voltage magnitude, line loss constraints, line flow constraints and
certain bus power constraints, the injection region and its convex hull have the same Pareto-front. Precisely,
the condition on the bus power constraints is: each bus is allowed
to have real and reactive power upper bounds, but two connected buses cannot both
simultaneously have real power lower bounds and there are no reactive power lower bounds. Through simulations with practical distribution networks, we show that these requirements are not stringent in actual operations. Independent works
\cite{Lavaei11c,Bose11} considered the OPF problem for a tree network,
although the authors there used the notion of load over-satisfaction
and did not consider thermal loss constraints.

The paper is organized as follows. In Section \ref{sec:II} we establish the notations, Section \ref{sec:nobounds} contains the result about the network with no operation constraints, Section \ref{sec:tree} contains theoretical and simulations results concerning trees, and Section \ref{sec:V} concludes the paper. The Appendices contain the results about non-tree networks and some of the proofs. 
%Due to space constraints, some proofs and details are omitted. A long version is posted on ArXiv \cite{Zhang11}.

The appendix address network with cycles. In some distribution systems, the network consists of a ring (cycle) feeder and tree networks hanging off the ring, therefore it is useful to understand the injection region of cycles. Ideally, one would like to state an analogous result as in
the tree network case. However, we could not yet prove such a strong
result. Instead, we characterize the convex hull of the voltage
magnitude constrained injection region if the network is a {\it
cycle with lossless links} and certain combinations of these networks with trees.

\section{Model and Notations} \label{sec:II}
We consider the AC power flow model so in general all variables are complex. Following the convention in power engineering, scalars representing voltage, current and power are denoted with capital letters. We use $\x$ to denote vectors, and $\bd{X}$ to denote matrices. $\x \odot \y$ denote the element-wise product between $\x$ and $\y$. Given two real vectors $\bd{x}$ and $\bd{y}$ of the same dimension, the notation $\bd{x} \leq \bd{y}$
denotes component-wise inequality and $\bd{x} < \bd{y}$ denotes
component-wise inequality with strict inequality in at least one
component. We denote Hermitian transpose by $(\cdot)^H$ and complex conjugation by $\conj(\cdot)$. We write $\bd{X} \sdp 0$ to mean $\bd{X}$ positive semidefinite.  Given a set $\mathcal{A} \subset \R^n$,
$\convhull(\mathcal{A})$ denote the convex hull of $\mathcal{A}$, i.e. the smallest convex set containing $\mc{A}$. 

Consider an electric network with $n$ buses. Throughout we assume the network is connected. We write $i \sim k$ if bus $i$ is connected to $k$, and $i \nsim k$ if they are not connected. Let $z_{ik}$ denote the complex impedance of the transmission line between bus $i$ and bus $k$, and $y_{ik}=\frac{1}{z_{ik}}=g_{ik}+j b_{ik}$. We have $g_{ik} >0$, and we assume that the lines are inductive (as in the Pi model) so $b_{ik} <0$. Note that $z_{ik}=z_{ki}$ and $y_{ik}=y_{ki}$. Let $z_{ii}$ ($y_{ii}$) denote the shunt impedance (admittance) of bus $i$ to ground. These shunt impedances can come from the capacitance to ground in the Pi model of the transmission line, the capacitor banks installed for reactive power injection, or modeling constant impedance loads. 
%Let $z_{ii}$ be the shunt element between node $i$ and ground, and $y_{ii}=\frac{1}{z_{ii}}$.
The bus admittance matrix is denoted by $\bd{Y}$ and defined as
\begin{equation} \label{eqn:Y}
Y_{ik}=
\begin{cases}
\sum_{l \sim i} y_{il}+ y_{ii} & \mbox{ if } i=k \\
-y_{ik} & \mbox { if } i \sim k \\
0 & \mbox { if } i \nsim k
\end{cases}.
\end{equation}
$\bd{Y}$ is symmetric. If the entries of $\bd{Y}$ are real, we say the network is purely resistive and if the entries are imaginary, we say the network is lossless. Lines in the transmission network are mainly inductive so it is sometimes assumed that the network is lossless. Let $\bd{v}=(V_1,V_2,\dots,V_n) \in \C^n$ be the vector of bus voltages and $\bd{i}=(I_1,I_2,\dots,I_n) \in \C^n$ be the vector of currents, where $I_i$ is the total current flowing out of bus $i$ to the rest of the network. By Ohm's law and Kirchoff's Current Law, $\bd{i}=\bd{Y}\bd{v}$. The complex power injected at bus $i$ is $S_i=P_i+jQ_i=V_i I_i^H$ where $P_i$ is the real power and $Q_i$ is the reactive power. A  positive $P_i$ means bus $i$ is generating real power and a negative $P_i$ means bus $i$ is consuming real power; similarly for $Q_i$. Let $\bd{p}=(P_1,P_2,\dots,P_n)$ be the vector of real powers and $\bd{q}=(Q_1,Q_2,\dots,Q_n)$ be the vector of reactive powers. 

The real power vector $\bd{p}=\Real(\bd{v} \odot \conj(\bd{i}))=\Real(\bd{v} \odot (\bd{Y}^H \bd{v}^H))=\Real(\diag(\bd{v} \bd{v}^H \bd{Y}^H))$ where $\diag(\bd{M})$ is the vector of diagonal elements of a matrix $\bd{M}$. Similarly, the reactive power vector $\bd{q}=\Imag (\diag(\bd{v}\bd{v}^H \bd{Y}^H))$. The resistive loss on a transmission line between buses $i$ and bus $k$ is given by $L_{ik}= |V_i -V_k|^2 g_{ik}$. The powers flowing from bus $i$ to bus $k$ is denoted $P_{ik}$ and $Q_{ik}$, and defined as
%\begin{equation} \label{eqn:flow}
$P_{ik}+j Q_{ik} =V_i |V_i-V_k|^* y_{ik}^*$.
%\end{equation}
Note $L_{ik}=P_{ik}+P_{ki}$.
\subsection{OPF Problem}
In power networks, we are often interested in solving the following OPF problem
\begin{subequations}
\label{eqn:opt}
\begin{align}
\mbox{minimize } & f(P_1,P_2,\dots,P_n)\\
\mbox{subject to } & \ul{V}_i \leq |V_i| \leq \ov{V}_i \label{eqn:optV} \\
                   & L_{ik} \leq l_{ik} \label{eqn:optL} \\
                   & P_{ik} \leq \ov{P}_{ik} \label{eqn:optPik} \\
                   & \ul{P}_i \leq P_i \leq \ov{P}_i \label{eqn:optPi} \\
                   & \ul{Q}_i \leq Q_i \leq \ov{Q}_i \label{eqn:optQi} \\ 
                   & \bd{p}+j \bd{q}=\Real(\bd{v} \bd{v}^H \bd{Y}^H), \label{eqn:optY}
\end{align}
\end{subequations}
where $f(P_1,P_2,\dots,P_n)$ is the cost function (not necessarily quadratic) defined on the real powers; \eqref{eqn:optV}, \eqref{eqn:optL}, \eqref{eqn:optPik}, \eqref{eqn:optPi} and \eqref{eqn:optQi} are the constraints corresponding to bus voltage, line thermal loss, line power flow and bus real and reactive power respectively; and \eqref{eqn:optY} is the physical law coupling voltage to power. The thermal loss constraints in \eqref{eqn:optL} are calculated from current rating of transmission lines and are usually the dominant constraints in distribution networks \cite{Kersting06}. Typically the data sheet of a line would have a maximum current rating $I_{\max}$ of the line, and this gives $l_{ik}=I_{\max}^2 R$, the maximum loss that can be tolerated across a line. In practice, $f$ is usually an increasing function of the power injections. For example, if $f(P_1,\dots,P_n)=P_1 + \cdots+ P_n$, then we are minimizing the loss in the network; or if $f$ is quadratic with positive coefficients, then we are minimizing the cost of generation.

In the rest of the paper we look at the feasible injection region, $\mathcal{P}$, defined as
\begin{align}
\mathcal{P} = \{ & \bd{p} \in \R^n : \bd{p}=\Real(\diag(\bd{v}\bd{v}^H \bd{Y}^H)), \ul{V}_i \leq |V|_i \leq \ov{V}_i \; \forall i, \nn \\ & L_{ik} \leq l_{ik} \; \forall i\sim k, P_{ik} \leq \ov{P}_{ik} \; \forall i \sim k, \nn \\ 
& \ul{P}_i \leq P_i \leq \ov{P}_i \forall i, \ul{Q}_i \leq Q_i \leq \ov{Q}_i\}.  \label{eqn:regionVPikbus}
\end{align}
Therefore $\mathcal{P}$ is the feasibility region of \eqref{eqn:opt}. Note the reactive powers are represented as a constraint of the injection region. This is because in most practical settings, the objective function of the optimization problem is in terms of {\em real} powers only. For example, the cost curve for an generator only includes the real power output; also, the consumers are only charged based on the amount of real power they consume (watt-hours). Since the objective function is in terms of real powers only, the injection region is the set of all real injections. 

\section{Network with No Operation Constraints} \label{sec:nobounds}
To warm up, let us first consider a network with no operation
constraints. Since there are no constraints, the injection region is
defined as
\begin{equation}
\mathcal{P}=\{\bd{p} \in \R^n : \bd{p}=\Real(\diag(\bd{v}\bd{v}^H \bd{Y}^H)).
\end{equation}
The reactive powers are ignored since we model reactive power as constraints in \eqref{eqn:opt}. In this case, the injection region has a simple characterization.
\begin{thm} \label{thm:1}
If the network is lossy\footnote{Every line has non-zero resistance}, then $\mathcal{P}$ is given by
\begin{equation}
\label{eq:energy}
\mathcal{P}= \{\bd{p} \in \R^n : \sum_{i=1}^n P_i > 0\} \cup \{\bd{0}\}.
\end{equation}
Therefore $\mathcal{P}$ is the union of the open upper half space of $\R^n$ and the origin $\bd{0}$. Note this region is connected and convex. If the network is lossless, then $\mathcal{P}$ is given by
\begin{equation}
\mathcal{P}=\{\bd{p} \in \R^n : \sum_{i=1}^n P_i=0\}.
\end{equation}
Therefore $\mathcal{P}$ is a hyperplane through the origin.
\end{thm}
This result is intuitive pleasing since it says if there are no
constraints in the network then the injection region is only limited
by the law of conservation of energy. Conservation of energy gives
the bound $\sum_{i=1}^n P_i \geq 0$, and if the network is not
lossless then $\sum_{i=1}^n P_i > 0$ except when all voltages are equal. In this case, all injections are $0$ so $\bd{p}=\bd{0}$. Theorem \ref{thm:1} states this is the only constraint on the injection region. The authors in \cite{Jarjis81, Hogan92} conjectured that the unconstrained injection region is convex, and \eqref{eq:energy} shows this is indeed the case. To proof this theorem, it is necessary to show that for every vector $\bd{p} \in \mc{P}$, there exists a voltage $\bd{v}$ that achieves $\bd{p}$. The details are given in the Appendix.

In practice, some of the constraints in \eqref{eqn:opt} would be binding. For example, the voltages magnitudes at each bus are
bounded. Figure \ref{fig:2busV} shows the injection region of a two
bus network with fixed voltage magnitudes. The region is an ellipse
(without the interior).  Even in this simple case, we see that the
injection region is no longer convex. The next section is devoted to
the study of the effect of constraints on the injection regions of
tree networks and their implications to optimization problems.

\section{Tree Networks} \label{sec:tree}
\subsection{Pareto-Front of Injection Region}
In this section we consider the full problem in \eqref{eqn:opt} for a tree network. 
%Suppose $f$ is convex. Traditionally, \eqref{eqn:opt} is said to be convex if $\mathcal{P}$ is a convex set. However, the convexity of $\mathcal{P}$ is not the relevant geometric feature to consider there.
The relevant geometric objects are the Pareto-optimal
points of $\mathcal{P}$ defined as:
\begin{defn}
Let $\mathcal{A} \subset \R^n$. A point $\bd{x} \in \mathcal{A}$ is
said to be a Pareto-optimal point if there does not exist another
point $\tilde{\bd{x}} \in \mathcal{A}$ such that $\tilde{\bd{x}} <
\bd{x}$. Denote the set of Pareto-optimal points of $\mathcal{A}$ as
$\mathcal{O}(A)$ and is sometimes called the Pareto-front of
$\mathcal{A}$.\footnote{Here we actually consider only the non-degenerative Pareto-optimal points. For a precise definition see \cite{Helton97}. In almost all applications, the set of degenerative Pareto-optimal points are of measure 0 and does not correspond to the minima of strictly increasing functions.}
\end{defn}
The Pareto-optimal points of $\mc{P}$ are of interest because only they can be the optimal solutions to \eqref{eqn:opt} when $f$ is increasing.
%\begin{defn} \label{def:Pconv}
%Let $\mathcal{A} \subset \R^n$. We say that the Pareto-front of $A$ is convex if $\mc{O}(\mc{A})=\mc{O}(\convhull(\mc{A}))$.
%\end{defn}
%Note this does not mean that $\mc{O}(\mc{A})$ as a set is convex, it means that the Pareto-front is invariant to taking the convex hull of $\mc{A}$. 
Under many circumstances, the Pareto-front of the injection region $\mc{P}$ is the same as the Pareto-front of $\convhull(\mc{P})$. Therefore, \eqref{eqn:opt}
is a convex optimization problem if $f$ is convex and increasing,
since we may replace the non-convex region $\mathcal{P}$ by a convex
region $\convhull(\mathcal{P})$ and obtain the same solutions.
Before stating the general result about the Pareto-front of $\mc{P}$
in Theorem \ref{thm:treelinebus}, it is instructive to use a two bus
example to see what are the Pareto-optimal points and the effect of
various kinds of constraints on them.

\begin{figure}[ht]%
\centering
\psfrag{V1}{$V_1$}
\psfrag{V2}{$V_2$}
\psfrag{P1}{$P_1$}
\psfrag{P2}{$P_2$}
\psfrag{y}{$y$}
\includegraphics[width=7cm]{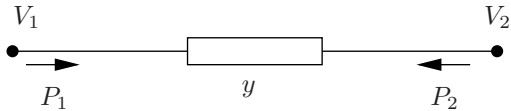}%
\caption{Two bus network.}%
\label{fig:2bus}%
\end{figure}

\begin{figure}[ht]
\centering
\subfigure[Voltage constrained.]{
\includegraphics[width=3.5cm]{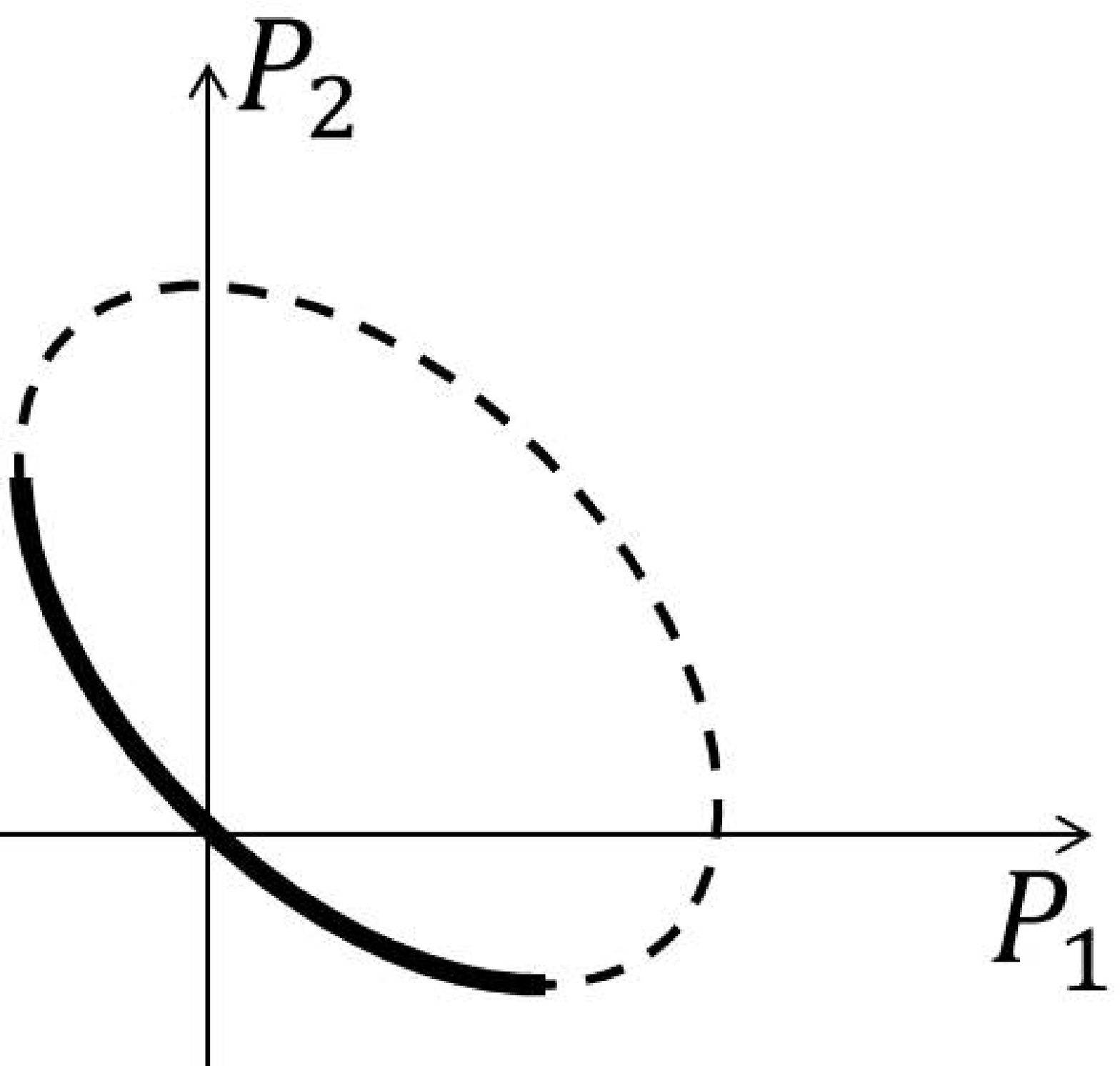}%
\label{fig:2busV}}
\subfigure[Voltage and loss constrained.]{
\includegraphics[width=3.5cm]{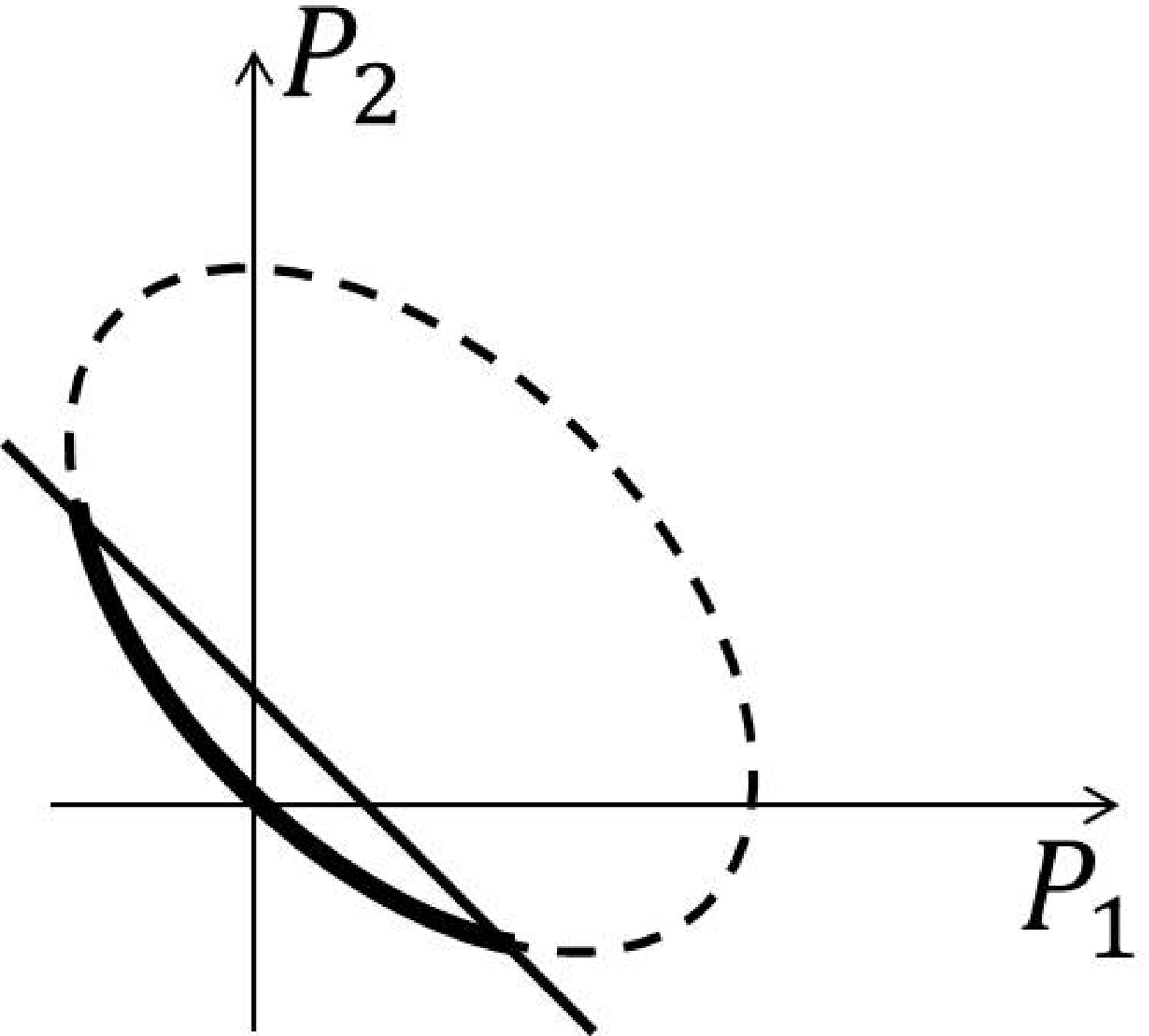}
\label{fig:2busVL}}
\caption{Voltage constrained and loss constrained injection regions. The parameters are $|V|_1=|V|_2=1, g=1,b=3$, all per unit.}
\label{fig:2busVVL}
\end{figure}

\begin{figure}[ht]
\centering
\subfigure[Both buses have power upper bounds.]{
\includegraphics[width=3cm]{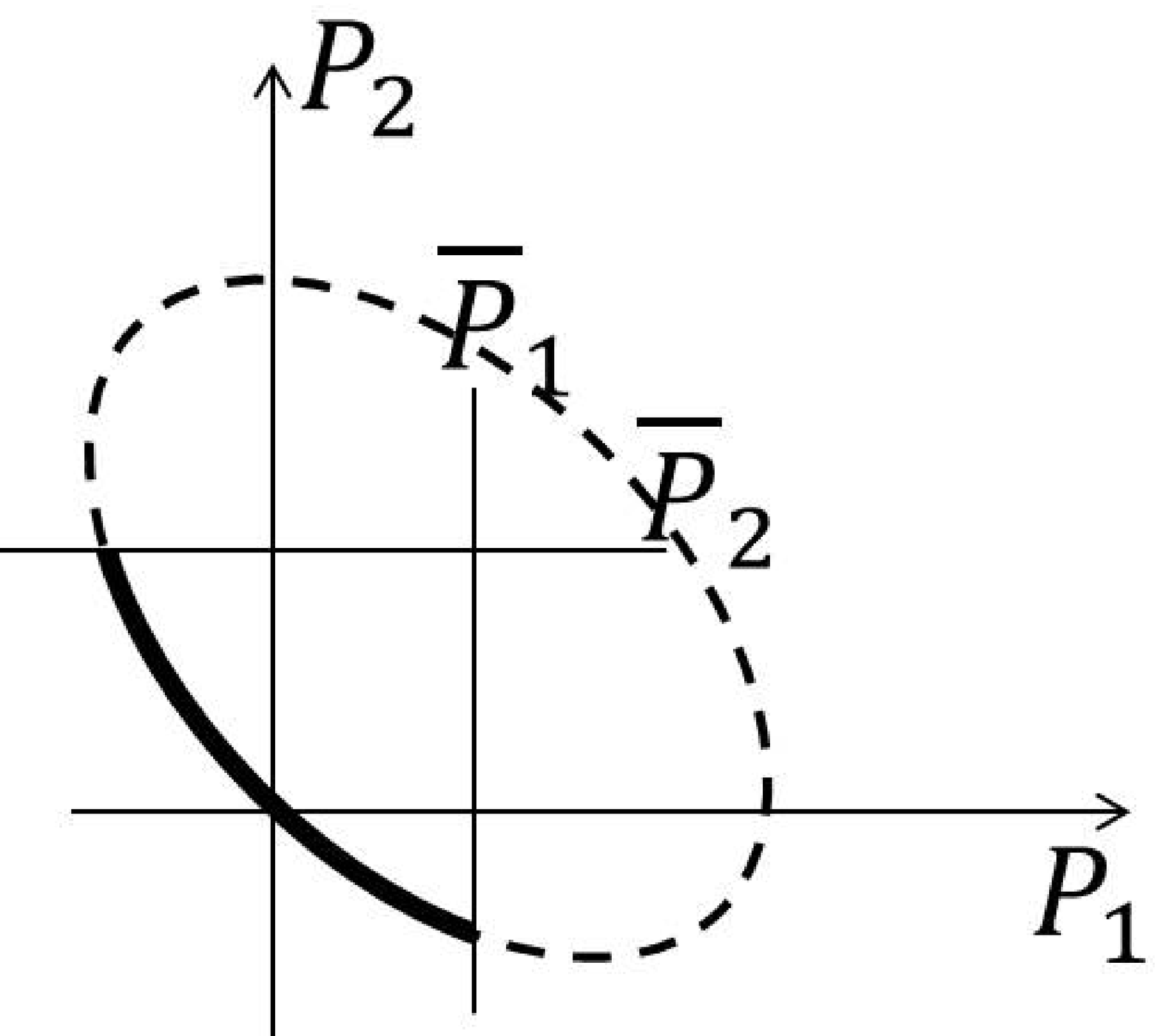}%
\label{fig:2busP1}}
\subfigure[$P_2$ have both upper and lower bounds.]{
\includegraphics[width=3cm]{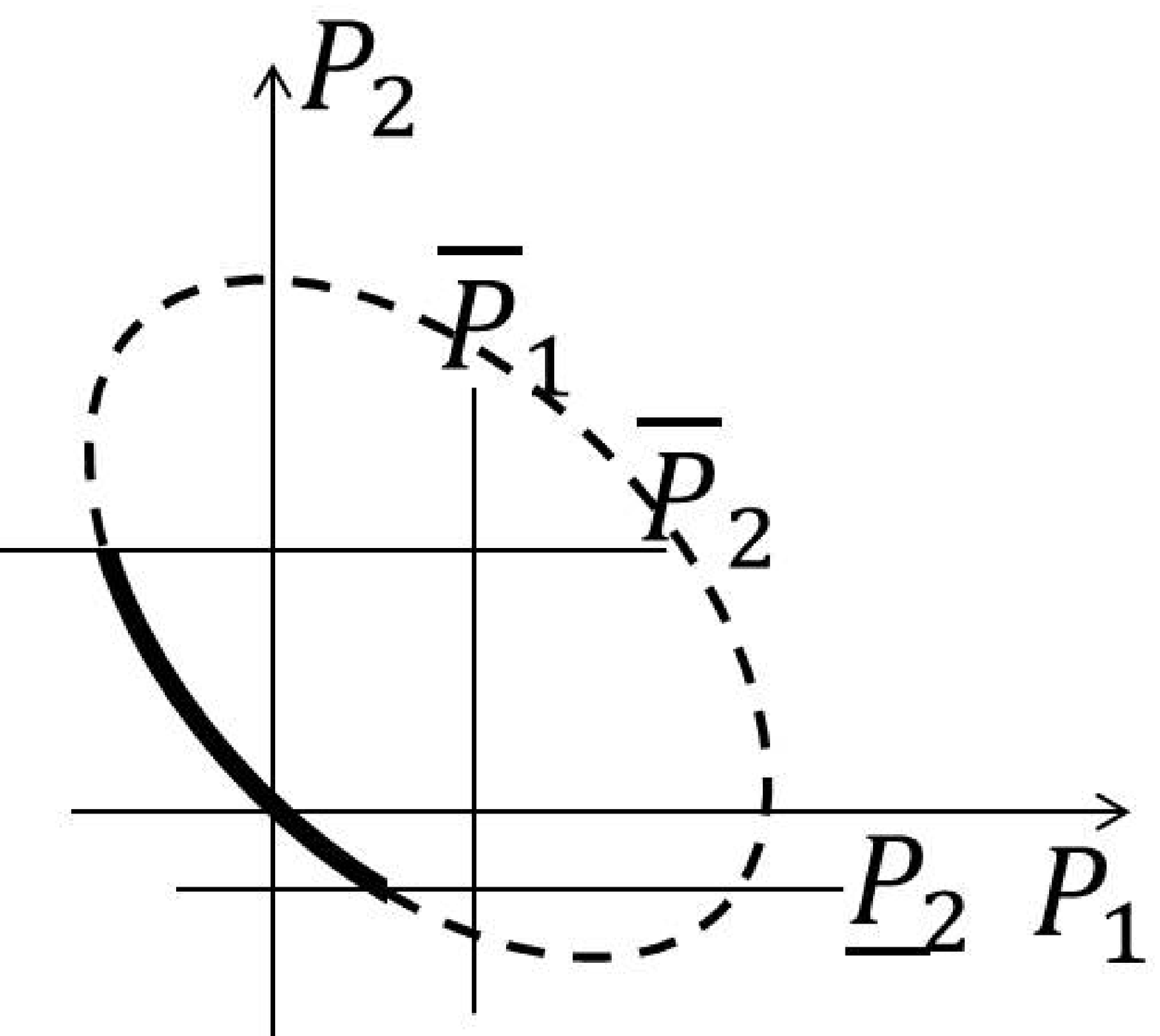}
\label{fig:2busP2}}
\subfigure[Both are lower bounded.]{
\includegraphics[width=3cm]{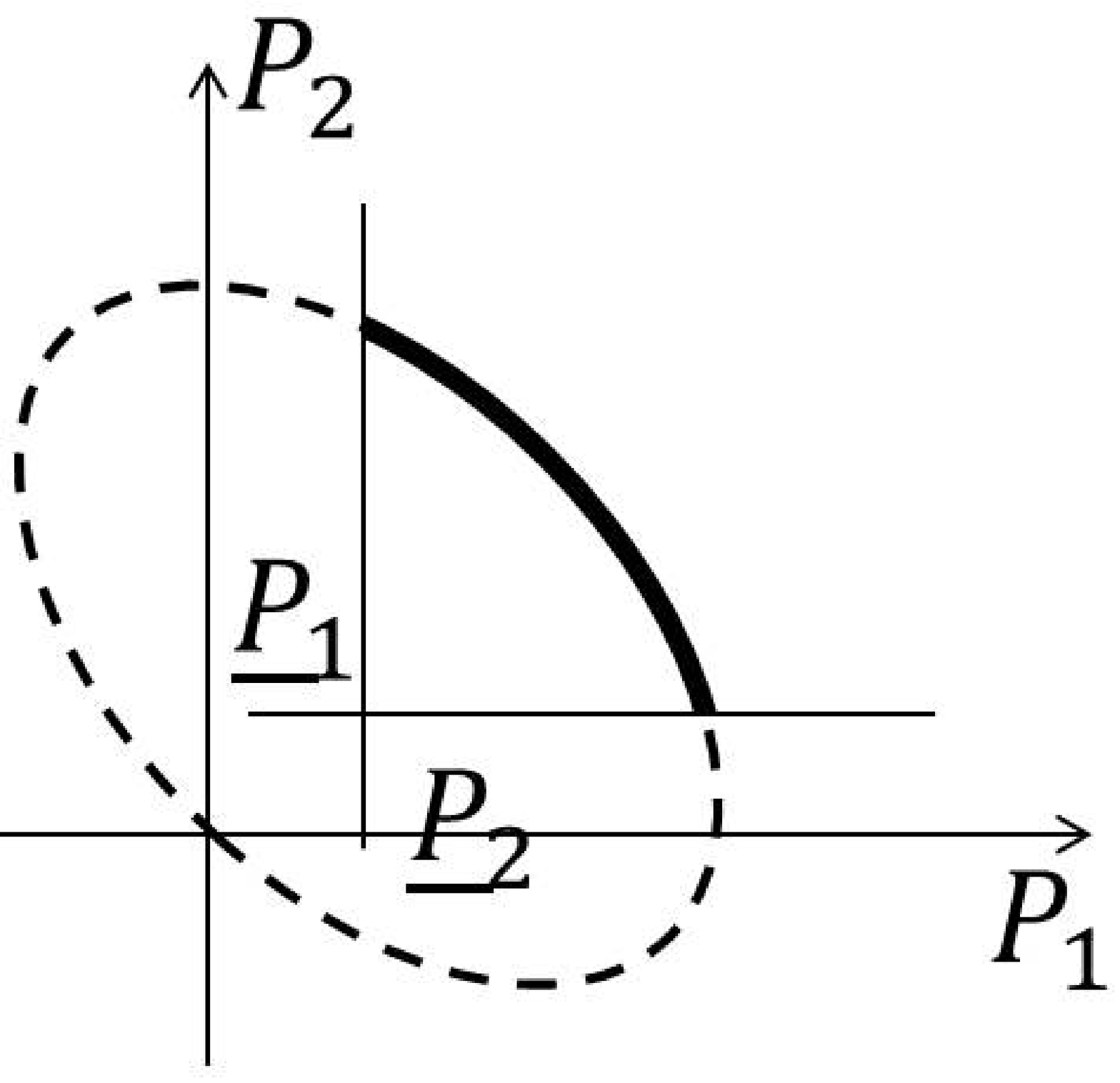}
\label{fig:2busP3}}
\caption{Three possible cases of the bus power constrained injection region.}
\label{fig:2busP}
\end{figure}

Consider the two bus example in Figure \ref{fig:2bus} where $y$ is the line admittance.
First consider the case where there are only voltage constraints.
Suppose that $|V_1|=|V_2|=1$ per unit. Then $\mathcal{P}$ is an
ellipse as shown in Figure \ref{fig:2busV}. The bold curve
represents the Pareto-front. Note $\convhull(\mc{P})$ is the filled ellipse. We can see that the Pareto-fronts of the empty and the filled ellipses are the same. Therefore, if we replace the non-convex empty ellipse by the convex filled ellipse in an optimization problem with increasing objective function, we would obtain the same solution. Next, we consider both voltage constraints and the loss constraint $P_{12}+P_{21}=P_1+P_2
\leq l$ for some $l$. This is presented by intersecting the ellipse
by a half plane as in Figure \ref{fig:2busVL}, and the bold curve is
the resulting Pareto-front, and we see that it is again the same as the Pareto-front of $\convhull(\mc{P})$. 
Next, consider both voltage and bus power constraints. In this case,
there are several possibilities, as represented in Figures
\ref{fig:2busP1}, \ref{fig:2busP2} and \ref{fig:2busP3}. In Figure
\ref{fig:2busP1}, both bus have power upper bounds, and the Pareto-front of $\mc{P}$ is the same as the Pareto-front of $\convhull(\mc{P})$. In Figure \ref{fig:2busP2}, $P_1$ has upper bound,
$P_2$ has both upper and lower bounds, and the Pareto-front of $\mc{P}$ is the same as the Pareto-front of $\convhull(\mc{P})$.  In Figure \ref{fig:2busP3}, both buses have lower bounds,
and we see that the Pareto-front of $\mc{P}$ is {\em not} the same as the Pareto-front of $\convhull(\mc{P})$.
Note that in the two bus case, the line flow constraints in
\eqref{eqn:optPik} correspond to Figure \ref{fig:2busP1}.

\begin{figure}[ht]
\centering
\subfigure[Real injection region.]{
\includegraphics[width=3cm]{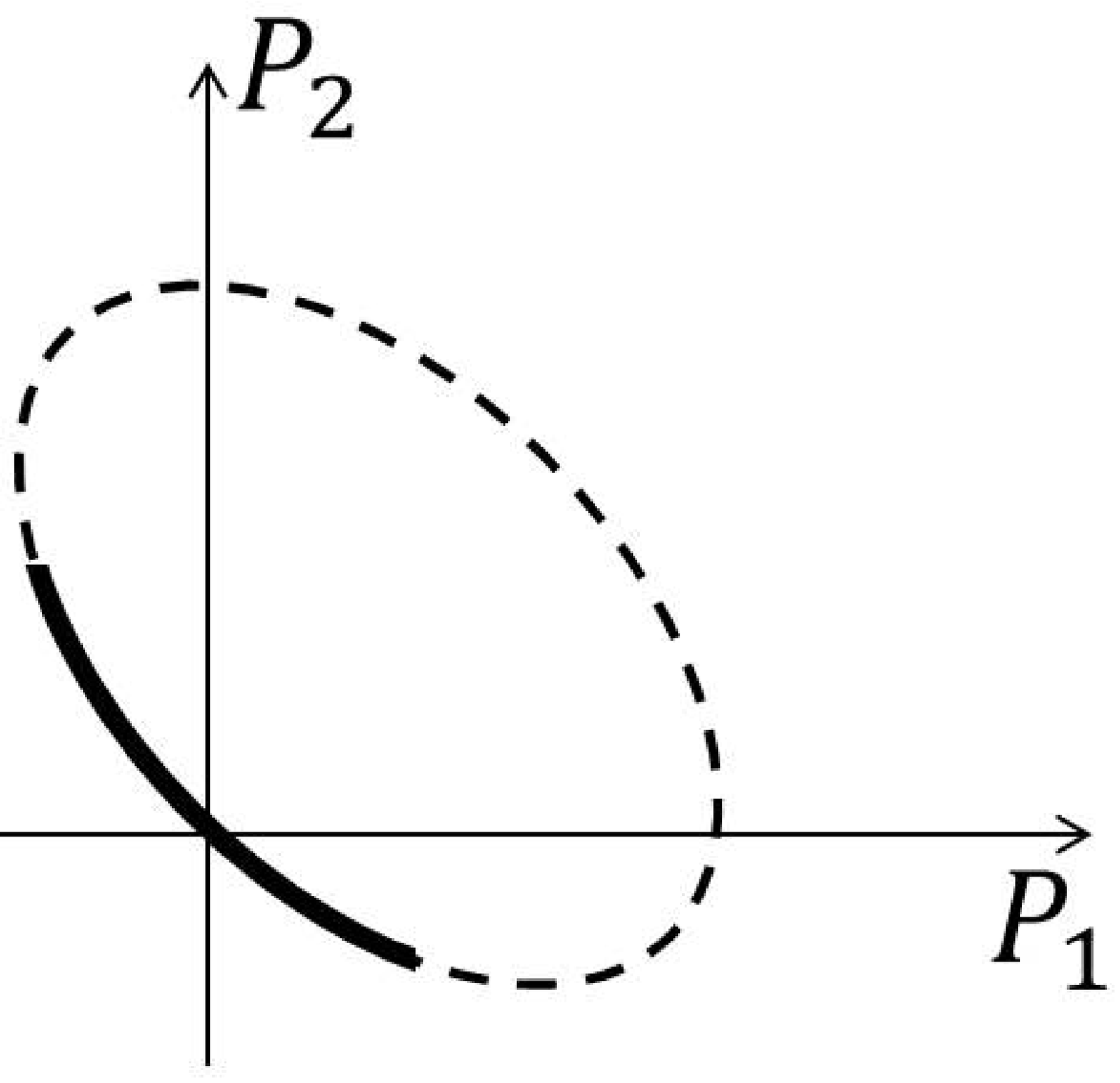}
\label{fig:2busPQ1}}
\subfigure[Reactive injection region.]{
\includegraphics[width=3cm]{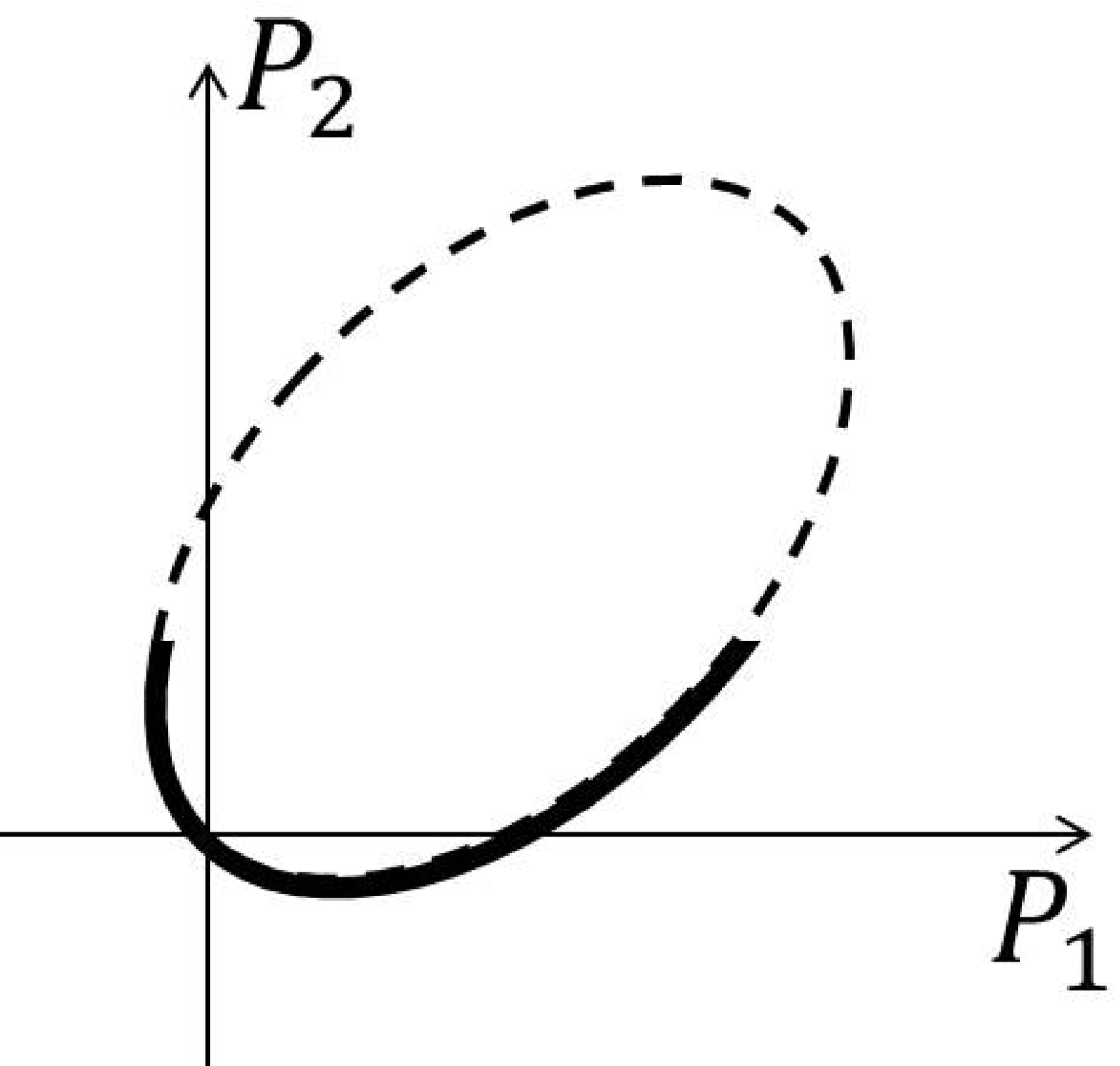}%
\label{fig:2busQ1}}
\subfigure[Real injection region.]{
\includegraphics[width=3cm]{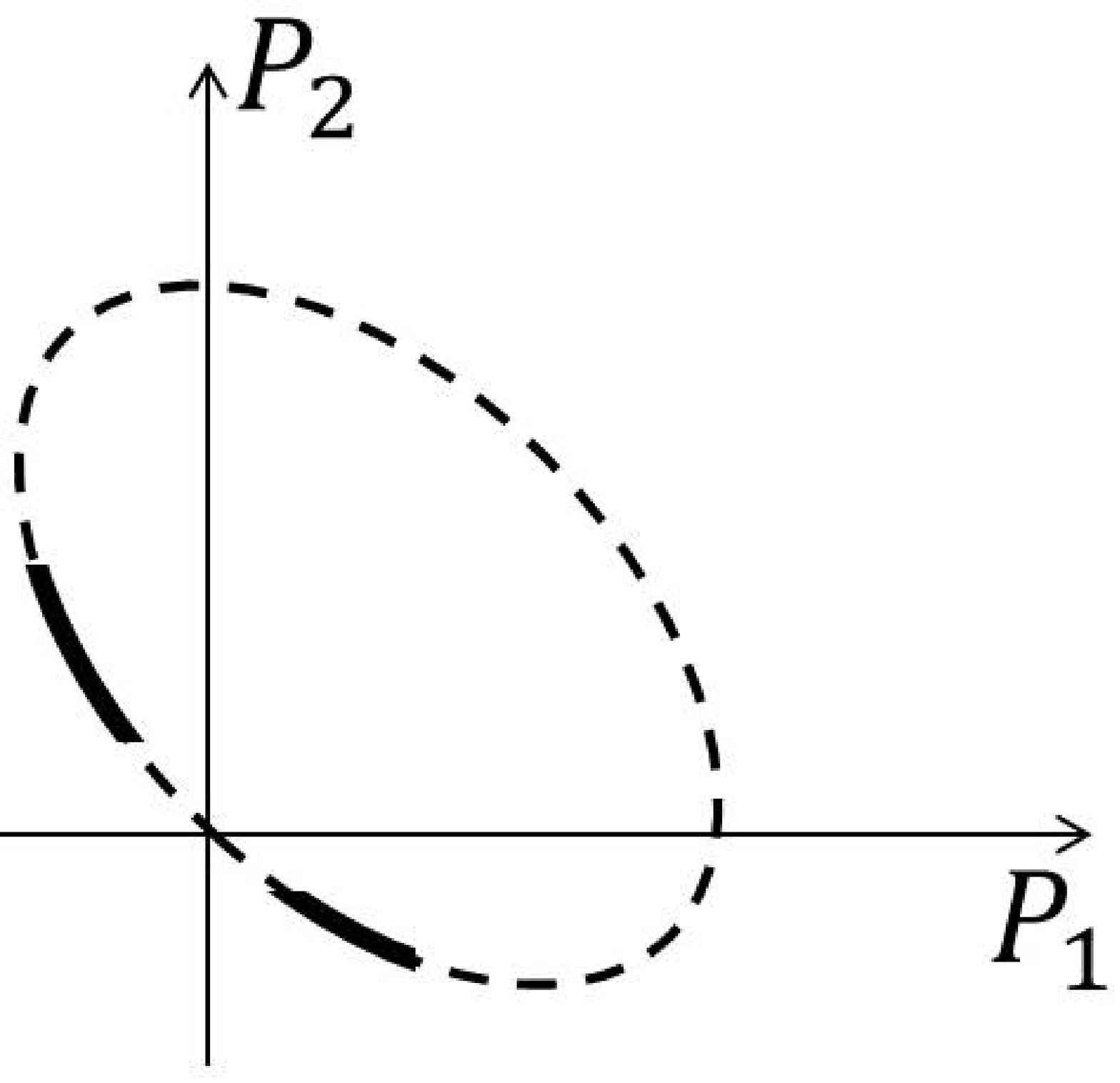}
\label{fig:2busPQ}}
\subfigure[Reactive injection region.]{
\includegraphics[width=3cm]{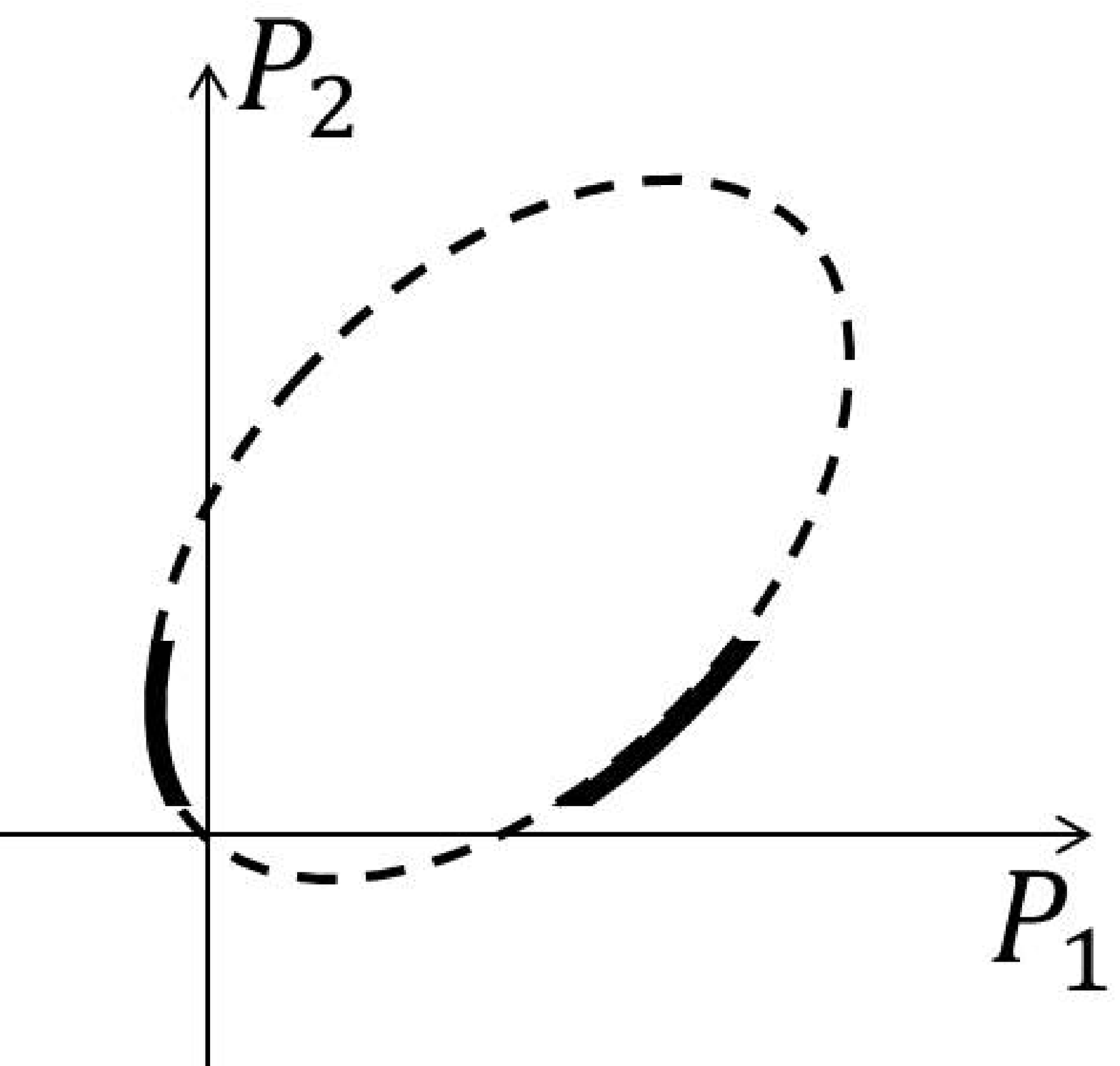}%
\label{fig:2busQ}}
\caption{Impact of reactive power constraints.}
\label{fig:2busPQQ}
\end{figure}

Next let us consider the effect of reactive power bounds. Figure \ref{fig:2busQ1} shows the feasible reactive power that can be achieved under the voltage constraint and the bold segment that satisfies the reactive power constraint $Q_2 \leq \ov{Q}_2$.  The bold segments in Figure \ref{fig:2busPQ1} shows the corresponding injection region. As we can see, the Pareto-front of $\mc{P}$ is the same as the Pareto-front of $\convhull(\mc{P})$. Next, Figure \ref{fig:2busQ} shows the bold segments that satisfies the constraint $\ul{Q}_2 \leq Q_2 \leq \ov{Q}_2$. As we can see, the Pareto-front of the Pareto-front of $\mc{P}$ is {\em not} the same as the Pareto-front of $\convhull(\mc{P})$ Therefore, in general we cannot extend the result to include reactive power lower bounds.

The intuition gained from the two bus example carries over for general
trees, and the general statement is given in Theorem
\ref{thm:treelinebus}.
\begin{thm} \label{thm:treelinebus}
Consider a tree network with $n$ buses.  Let the injection region
$\mathcal{P}$ defined as in \eqref{eqn:regionVPikbus}. Suppose two conditions are satisfied:
\begin{enumerate}
	\item If $i \sim k$, then either $\ul{P}_i = -\infty$ or $\ul{P}_k = -\infty$.
	\item $\ul{Q}_i=-\infty$ for all $i$. 
\end{enumerate}
%Under the condition if $i
%\sim k$ in the network, either the real power lower bound on bus $i$ is not tight or the real power bound on bus $k$ is not tight and suppose the reactive power lower bounds are all not tight. 
The Pareto-front of $\mc{P}$ is the same as the Pareto-front of $\convhull(\mc{P})$. 
\end{thm}
The condition on the bus power lower bounds means that if two buses are connected, then not both can have a tight bus real power lower bound. Also, the theorem requires that all the reactive lower bounds to be not tight. This can be seen as a generalization of the well known load over-satisfaction concept \cite{Baldick06}. In load over-satisfaction, all the lower bounds on real and reactive power are removed. But Theorem \ref{thm:treelinebus} states it is not necessary to remove all the lower bounds. 
%For example, consider the tree network in Figure \ref{fig:6tree}. Theorem \ref{thm:treelinebus} applies if the following sets of buses have tight real power lower bounds: $\{1,4,5,6\},\{2,3\},\{2,6\},\{3,4,5\}$.
%\begin{figure}[ht]%
%\centering
%\includegraphics[height=2cm]{6tree.eps}%
%\caption{A 3 level tree with $6$ buses.}%
%\label{fig:6tree}%
%\end{figure}

\begin{IEEEproof}
To prove the theorem, first we define an optimization problem in term of the injection region. In this optimization problem, we want to write every quantity as a quadratic form of the complex voltages.

The resistive loss on the transmission line between buses $i$ and $k$ can be written as $L_{ik}= \bd{v}^H \bd{G}_{ik} \bd{v}$ where $\bd{G}_{ik}$ is a matrix with the $(i,i)$th entry and the $(k,k)$th entry being $g_{ik}$, and the $(i,k)$th entry and the $(k,i)$th entry being $-g_{ik}$ and all other entries being $0$. The power flow from bus $i$ to bus $k$ can be written as $P_{ik}=\bd{v}^H \bd{A}_{ik} \bd{v}$, where $\bd{A}_{ik}$ is a matrix with $(i,i)$th entry $g_{ik}$, the $(i,k)$th entry $\frac12 (-g_{ik} -j b_{ik})$, the $(k,i)$th entry $\frac12 (-g_{ik}+jb_{ik})$ and all the other entries $0$. Let $\bd{A}_i= \frac12 (\bd{E}_i \bd{Y}+\bd{Y}^H \bd{E}_i)$ where $\bd{E}_i$ is the diagonal matrix with $1$ at the $(i,i)th$ entry and $0$ everywhere else. Similarly let $\bd{B}_i=\frac{1}{2j} (\bd{Y}^H \bd{E}_i -\bd{E}_i \bd{Y})$. Then the powers injected at bus $i$ is given by $P_i=\bd{v}^H \bd{A}_i \bd{v}$ and $Q_i=\bd{v}^H \bd{B}_i \bd{v}$.

Consider the following optimization problem
\begin{align}\label{eqn:opfPVlinebus}
J=\mbox{minimize } & \sum_{i=1}^n c_i P_i\\
 \mbox{subject to } & \ul{V}_i \leq |V_i| \leq \ov{V}_i, \; \forall i \nn \\
& \bd{v}^H \bd{G}_{ik} \bd{v} \leq l_{ik} \; \forall i \sim k \nn \\
& \bd{v}^H \bd{A}_{ik} \bd{v} \leq \ov{P}_{ik} \; \forall i \sim k \nn \\
& \ul{P}_i \leq \bd{v}^H \bd{A}_i \bd{v} \leq \ov{P}_i \nn \\
& \ul{Q}_i \leq \bd{v}^H \bd{B}_i \bd{v} \leq \ov{Q}_i \nn \\
&  \bd{p}+j \bd{q} =\diag(\bd{v}\bd{v}^H \bd{Y}^H). \nn
\end{align}

The $c_i$'s can be interpreted as the costs of the power generation and \eqref{eqn:opfPVlinebus} is an optimal power flow problem with a linear cost function. To expose the potential non-convexity, we can equivalently write it as
\begin{align} \label{eqn:opfWlinebus}
 J=\mbox{minimize } & \sum_{i=1}^n c_i P_i \\
 \mbox{subject to } & \ul{V}_i^2 \leq W_{ii} \leq \ov{V}_i^2, \; \forall i \nn \\
 & \Tr(\bd{G}_{ik} \bd{W}) \leq l_{ik} \; \forall i \sim k \nn \\
 & \Tr(\bd{A}_{ik} \bd{W}) \leq \ov{P}_{ik} \; \forall i \sim k \nn \\
 & \ul{P}_i \leq \Tr(\bd{A}_i \bd{W}) \leq \ov{P}_{i} \nn \\
 & \ul{Q}_i \leq \Tr(\bd{B}_i \bd{W}) \leq \ov{Q}_i \nn \\
 & \bd{p}+j \bd{q} =\diag(\bd{W} \bd{Y}^H) \nn \\
 & \bd{W} \sdp 0 \nn \\
 & \rank(\bd{W})=1, \nn
\end{align}
where $\bd{W}=\bd{v}\bd{v}^H$ and the non-convexity enters as the rank 1 constraint on $\bd{W}$. Relaxing this rank $1$ constraint and eliminating $\bd{p}$ and $\bd{q}$, we get
\begin{align} \label{eqn:opfrWlinebus}
J_1=\mbox{minimize } & \Tr(\bd{M} \bd{W}) \\
 \mbox{subject to } & \ul{V}_i^2 \leq W_{ii} \leq \ov{V}_i^2, \; \forall i \nn \\
 & \Tr(\bd{G}_{ik} \bd{W}) \leq l_{ik} \; \forall i \sim k \nn \\
 & \Tr(\bd{A}_{ik} \bd{W}) \leq \ov{P}_{ik} \; \forall i \sim k \nn \\
 & \ul{P}_i \leq \Tr(\bd{A}_i \bd{W}) \leq \ov{P}_{i}, \nn \\
 & \ul{Q}_i \leq \Tr(\bd{B}_i \bd{W}) \leq \ov{Q}_i, \nn \\
 & \bd{W} \sdp 0 \nn
\end{align}
where $\bd{M}=\frac12 (\bd{C}\bd{Y}+\bd{Y}^H\bd{C})$ and $\bd{C}=\diag(c_1,\dots,c_n)$. Note $M$ is Hermitian.

Geometrically, the relaxation from \eqref{eqn:regionVPikbus} to \eqref{eqn:opfrWlinebus} enlarges the feasible injection region to a convex region given by
\begin{align}
\tl{\mathcal{P}}=\{& \bd{p}: \bd{p}=\Real(\diag(\bd{W}\bd{Y}^H)), \ul{V}_i^2 \leq W_{ii} \leq \ov{V}_i^2 \; \forall i, \\
& \Tr(\bd{G}_{ik}\bd{W}) \leq l_{ik} \; \forall i \sim k, \Tr(\bd{A}_{ik} \bd{W}) \leq \ov{P}_{ik} \; \forall i \sim k, \nn \\
& \ul{P}_i \leq \Tr(\bd{A}_i \bd{W}) \leq \ov{P}_i, \ul{Q}_i \leq \Tr(\bd{B}_i \bd{W}) \leq \ov{Q}_i, \bd{W} \sdp 0\} \nn.
\end{align}
We want to show that the two regions have the same Pareto-front.
That is,
$\mathcal{O}(\mathcal{P})=\mathcal{O}(\tl{\mathcal{P}})$. Since
$\tl{\mathcal{P}}$ is convex, its Pareto-front is easily
explored. Note in general $\tl{\mathcal{P}} \supseteq
\convhull(\mathcal{P})$ and the inclusion can be strict. However, if
$\mathcal{P}$ and $\tl{\mathcal{P}}$ have the same Pareto-front,
then so does $\convhull(\mathcal{P})$. 

The proof of the theorem follows from the following claim.
\begin{claim} \label{clm:tree}
Suppose $c_i > 0$ for all $i$. Then the optimal solution to \eqref{eqn:opfrWlinebus} is unique and has rank $1$ if for every connected pair of buses $(i,k)$ in the network, one of them do not have tight bus power lower bound, and all reactive power lower bounds are not tight. 
\end{claim}
This claim is a stronger statement then saying $J=J_1$, it also
states that the optimal solution to the relaxed solution is unique.
Assuming for now the claim is true. Then since $\tl{\mathcal{P}}$ is
convex, we can explore its Pareto-front by linear functions with
positive costs \cite{Boyd04}. More precisely, a point
$\tl{\bd{p}} \in \tl{\mathcal{P}}$ is a Pareto-optimal if and
only if it is an optimal solution to \eqref{eqn:opfrWlinebus} for
some positive costs. From the claim, all the optimal solutions are
achieved by a $\bd{W}$ of rank $1$, therefore they can be achieved
by using a voltage vector $\bd{v}$. Therefore if $\bd{p} \in
\tl{\mathcal{P}}$ is a Pareto-optimal, then $\bd{p} \in
\mathcal{P}$. Since $\tl{\mathcal{P}} \supseteq \mathcal{P}$,
$\bd{p}$ is also a Pareto-optimal point of $\mathcal{P}$. So
$\mathcal{O}(\mathcal{P}) \supseteq
\mathcal{O}(\tl{\mathcal{P}})$. To show the other direction,
suppose there exists a point $\bd{p} \in \mathcal{O}(\mathcal{P})$
but not in $\mathcal{O}(\tl{\mathcal{P}})$. Then there is a
point $\tilde{\bd{p}} \in \mathcal{O}(\tl{\mathcal{P}})$ such
that $\tilde{\bd{p}} \leq \bd{p}$. But $\tilde{\bd{p}} \in
\mathcal{O}(\mathcal{P})$, contradicting the fact $\bd{p}$ is a
Pareto-optimal point of $\mathcal{P}$. Therefore
$\mathcal{O}(\mathcal{P}) \subseteq
\mathcal{O}(\tl{\mathcal{P}})$ and thus
$\mathcal{O}(\mathcal{P}) = \mathcal{O}(\tl{\mathcal{P}})$. It
remains to prove claim \ref{clm:tree}.

We are to show  that the optimal solution to \eqref{eqn:opfrWlinebus}, $\bd{W}^*$, is rank 1. We do this through duality theory. The dual of \eqref{eqn:opfrWlinebus} is
\begin{align}
\mbox{maximize } & \sum_{i=1}^n (\ul{\la}_i \ul{V}_i^2 -\ov{\la}_i \ov{V}_i^2)-\sum_{i \sim k} \mu_{ik} l_{ik} - \nn \\
 & \sum_{i \sim k} (\nu_{ik} \ov{P}_{ik}+\nu_{ki} \ov{P}_{ki}) +  \sum_{i=1}^n (\ul{\sigma}_i \ul{P}_i -\ov{\sigma}_i \ov{P}_i- \rho_i \ov{Q}_) \nn \\
 \mbox{subject to } & \bd{\La} + \sum_{i \sim k} \mu_{ik} \bd{G}_{ik}+ \sum_{i \sim k} (\nu_{ik} \bd{A}_{ik}  + \nu_{ki} \bd{A}_{ki}) \nn \\
 &+ \sum_{i=1} (\sigma_i \bd{A}_i + \rho_i \bd{B}_i) + \bd{M} \sdp 0, \label{eqn:opfDlinebus}
\end{align}
where $\ov{\la}_i$ and $\ul{\la}_i$ are the Lagrange multiplier associated with the voltage upper and lower bounds and $\la_i=\ov{\la}_i-\ul{\la}_i$ and $\bd{\La}=\diag(\la_1,\dots,\la_n)$, $\mu_{ik}$ are the Lagrange multiplier associated with the thermal constraints, $\nu_{ik}$ and $\nu_{ki}$ are the Lagrange multipliers associated with the flow constraints, and $\ov{\sigma}_i$ and $\ul{\sigma}_i$ are the Lagrange multiplier associated with the power upper and lower bounds and $\sigma=\ov{\sigma}_i-\ul{\sigma}_i$. Since we assume that the reactive power lower bounds constraints are not tight, $\rho_i$ is the Lagrange multiplier associated with the reactive power upper bounds. Note \eqref{eqn:opfDlinebus} is also the dual of \eqref{eqn:opfPVlinebus} so the gap between $J$ and $J_1$ is called the duality gap.

Let $\widetilde{\bd{M}}=\sum_{i \sim k} (\mu_{ik} \bd{G}_{ik})+\sum_{i \sim k} (\nu_{ik} \bd{A}_{ik}+ \nu_{ki} \bd{A}_{ki})+\sum_{i=1}^n (\sigma_i \bd{A}_i+ \rho_i \bd{B}_i) + \bd{M}$. Let $\bd{W}^*$ denote the optimal solution of \eqref{eqn:opfrWlinebus} and $\bd{\La}^*$ the optimal solution of \eqref{eqn:opfDlinebus}, by the complimentary slackness condition \cite{Boyd04},
\begin{equation}\label{eqn:cs}
\Tr((\bd{\La}^*+\widetilde{\bd{M}})\bd{W}^*) =0.
\end{equation}
Since both $\bd{W}^*$ and $\bd{\La}^*+\widetilde{\bd{M}}$ are positive semidefinite, \eqref{eqn:cs} implies that $(\bd{\La}^*+\widetilde{\bd{M}})\bd{W}^*=0$. Therefore $\bd{W}^*$ is in the null space of $\bd{\La}^*+\widetilde{\bd{M}}$ and $\rank(\bd{\La}^*+\widetilde{\bd{M}})+\rank(\bd{W}^*) \leq n$. So to show $\rank(\bd{W}^*)=1$ it suffices to show $\rank(\bd{\La}^*+\widetilde{\bd{M}}) \geq n-1$. This is done by considering the topology of the network and thus the structure of $\widetilde{\bd{M}}$.

Given a $n \times n$ matrix $\bd{A}$ and a graph $G$ with $n$ nodes, we say that $\bd{A}$ fits $G$ if for $i \neq k$, $A_{ik} = 0$ if and only if $(i,k)$ is not an edge in $G$. The values on the diagonal of $\bd{A}$ are unconstrained. The next lemma from \cite{Holst03} relates the topology of a graph and the rank of matrix that fits it.
\begin{lem}[Theorem 3.4 in \cite{Holst03}] \label{lem:treefit}
Let $G$ be a graph that is a connected tree of $n$ nodes. Suppose $\bd{A}$ is a $n \times n$ complex positive semidefinite matrix that fits $G$. Then $\rank(\bd{A}) \geq n-1$.
\end{lem}
We want to apply this lemma to the matrix $\bd{\La}^*+\widetilde{\bd{M}}$. Since $\bd{\La}^*$ is diagonal, only $\widetilde{\bd{M}}$ matters and its $(i,k)$th entry, $\widetilde{M}_{ik}$ is given by
\begin{equation*}
\begin{cases}
-\frac12 ((c_i+c_k+\mu_{ik}+ \nu_{ik}+\nu_{ki}+\sigma_i+\sigma_k) g_{ik} -\rho_i b_{ik} \\
 \mbox{       } + j (c_i-c_k+\mu_{ik}-\mu_{ki}+\sigma_i-\sigma_k) b_{ik}+\rho_i g_{ik}) \mbox{ if } i \sim k \\
0  \mbox{ if } i \nsim k
\end{cases}
\end{equation*}
Therefore if $\widetilde{M}_{ik}=0$ if bus $i$ is not connected to bus $k$. For $\widetilde{\bd{M}}$ to fit the network, $\widetilde{M}_{ik}$ needs to be nonzero if $i$ is connected to $k$.

If $i \sim k$, for $\widetilde{\bd{M}}_{ik}$ to be zero we need
\begin{align}
(c_i+c_k+\mu_{ik}+\nu_{ik}+\nu_{ki}+\sigma_i+\sigma_k)g_{ik}-\rho_i b_{ik} &= 0 \label{eqn:11} \\
(c_i-c_k+\nu_{ik}-\nu_{ki}+\sigma_i-\sigma_k)b_{ik} + \rho_i g_{ik} =0. \label{eqn:22}
\end{align}
Multiplying \eqref{eqn:11} by $g_{ik}$ and \eqref{eqn:22} by $b_{ik}$ and adding we get 
\begin{align*}
 (c_i+c_k+\mu_{ik}+\nu_{ik}+\nu_{ki}+\sigma_i+\sigma_k) &g_{ik}^2 \\
 +(c_i-c_k+\nu_{ik}-\nu_{ki}+\sigma_i-\sigma_k) &b_{ik}^2=0
\end{align*}
We are to show that $(c_i+c_k+\mu_{ik}+\nu_{ik}+\nu_{ki}+\sigma_i+\sigma_k)=(c_i-c_k+\nu_{ik}-\nu_{ki}+\sigma_i-\sigma_k)=0$. If this not the case, then suppose $(c_i+c_k+\mu_{ik}+\nu_{ik}+\nu_{ki}+\sigma_i+\sigma_k)<0$ and $(c_i-c_k+\nu_{ik}-\nu_{ki}+\sigma_i-\sigma_k)>0$. But $\rho_i \geq 0$ since it is a Lagrange multiplier and $g_{ik} > 0$, this contradicts \eqref{eqn:22}. Similarly, since $b_{ik} < 0$ (lines are inductive), we cannot have $(c_i+c_k+\mu_{ik}+\nu_{ik}+\nu_{ki}+\sigma_i+\sigma_k)>0$ and $(c_i-c_k+\nu_{ik}-\nu_{ki}+\sigma_i-\sigma_k)<0$. Therefore we get the simultaneous equations in \eqref{eqn:1} and \eqref{eqn:2}.   
\begin{align}
c_i+c_k+\mu_{ik}+\nu_{ik}+\nu_{ki}+\sigma_i+\sigma_k &= 0 \label{eqn:1} \\
c_i-c_k+\nu_{ik}-\nu_{ki}+\sigma_i-\sigma_k =0. \label{eqn:2}
\end{align}
Note $\mu_{ik}$, $\nu_{ik}$ and $\nu_{ki}$ are always nonnegative since they are the Lagrange multipliers associated with upper bounds. Suppose the bus power lower bound is not tight for bus $i$, then $\sigma_i \geq 0$. Adding \eqref{eqn:1} with \eqref{eqn:2} gives $2 c_i + \mu_{ik}+2 \nu_{ik}+2 \sigma_i=0$, this is not possible since $c_i >0$. On the other hand, suppose the bus power lower bound is not tight for bus $k$, then $\sigma_k \geq 0$. Subtracting \eqref{eqn:2} from \eqref{eqn:1} gives $2 c_k +\mu_{ik}+2 \nu_{ki}+2 \sigma_k =0$, which is not possible since $c_k >0$. Therefore $\widetilde{\bd{M}}$ fits a connected tree. Now apply Lemma \ref{lem:treefit} to the matrix $\bd{\La}^*+\widetilde{\bd{M}}$ gives $\rank(\bd{\La}^*+\widetilde{\bd{M}}) \geq n-1$, therefore $\rank(\bd{W}^*) \leq 1$. If the problem is feasible, then $\rank(\bd{W}^*)=1$.
\end{IEEEproof}

The authors in \cite{Lavaei11} showed that there is no gap if the network is purely resistive and all costs positive. Interpreting this in our language, they showed that the Pareto-front of the injection region of the resistive network is the same as that of its convex hull. In contrast, our results are based on the topology of the network, and do not need to make assumption that the network is purely resistive.

\subsection{Simulation Results}
In this section, we consider the voltage support problem in distribution networks. Due to the emergence of renewable generations and the high $R/X$ ratio in distribution networks, this is an interesting and non-trivial problem. Here we take the objective to be minimizing the total resistive loss in the network. So $f(P_1,\dots,P_n)=\sum_{i=1}^n P_i$, and the relaxed optimization problem in \eqref{eqn:opfrWlinebus} becomes
\begin{subequations}
\label{eqn:optVol}
\begin{align} 
J_1=\mbox{minimize } & \Tr(\bd{M} \bd{W}) \\
 \mbox{subject to } &  W_{ii} = \ov{V}_i^2 \; \forall i  \\
 & \Tr(\bd{G}_{ik} \bd{W}) \leq l_{ik} \; \forall i \sim k \label{eqn:thermalVol} \\
 & \ul{P}_i \leq \Tr(\bd{A}_i \bd{W}) \leq \ov{P}_{i}  \label{eqn:PVol} \\
 & \ul{Q}_i \leq \Tr(\bd{B}_i \bd{W}) \leq \ov{Q}_i \label{eqn:QVol} \\
 & \bd{W} \sdp 0,
\end{align}
\end{subequations}
where $M=\frac12 (Y^H+Y)$ and $\ov{V}_i$ is the given voltage level that we want to support. We obtain the test networks from the distribution network database in \cite{testfeeders}. In these test networks, the transmission line data and a typical power consumption profile is presented. From the transmission line data we obtain the $\bd{Y}$ matrix, and the thermal limits in \eqref{eqn:thermalVol} can be obtained from the maximum current ratings (line power flow rating was not included in the datasheets). We take $\ov{V}_i$ to be $1 p.u.$ for all buses. 

To verify our result, we need to construct the lower and upper bounds on $P_i$'s and $Q_i$'s. We assume the feeder acts like a slack bus, so it does not have any real or reactive power constraints. We consider two ways to construct the constraints for the other buses. One is that we assume a medium level penetration of solar generation at each bus. Let $\tl{P}_i$ be the typical real power consumption reported in \cite{testfeeders}, we randomly generate $\ov{P}_i \in [\tl{P}_i, 1.2 \tl{P}_i]$ and $\ul{P}_i \in [0.8 \tl{P}_i, \tl{P}_i]$. That is, we assume that the solar penetration level is about 20\% of the current power consumption, and depending on the environmental conditions, a real time $\ov{P}_i$ and $\ov{Q}_i$ is realized. Let $\tl{Q}_i$ be the typtical reactive power consumption of the network, we assume that $\ul{Q}_i=0$ and $\ov{Q}_i=1.2 \tl{Q}_i$. Note these bounds are typically fixed since they are provided by the power eletronics on the solar cells and is not dependent on the radiation levels. The newest power electronics availiable now have the ability to adjust its reactive power output within some bounds. We choose the lower bounds to be $0$ because all the power electronics can be adjusted to output $0$ reactive power. If a test case is generated this way, we say it is a nominal case since it came from a nominal operating point. If the parameters are choosen this way, all nodes except the feeder are withdrawing real power from the network. Rooting the tree at the feeder,  all real power flows in one direction: from the feeder to the leaf buses. 

Another way to generate the upper lower bounds is to randomly draw them such that $-2 \tl{P}_i < \ul{P}_i < \ov{P}_i < 2 \tl{P}_i$ and $-2 \tl{Q}_i < \ul{Q}_i < \ov{Q}_i < 2 \tl{Q}_i$. Note the problem parameter chosen this way may not correspond to any practical operation conditions. There could be multiply nodes with positive power injections into the network, resulting in real power flows that are bidirectional. We call this case the random case. 

During the simulations we solve the relaxed convex problem in \eqref{eqn:optVol}. We are interested in when the relaxed problem is tight; that is, when the optimal solution $\bd{W}^*$ to \eqref{eqn:optVol} is rank 1. We consider 3 networks, the 8-bus, 13-bus and the 34-bus networks. For the each of the networks, we run 1000 instances of the nominal and random generated cases. Table \ref{tab:sim} shows the number of times that $\bd{W}^*$ is rank 1 out of 1000 times. 
\begin{table}[ht]
\centering
\begin{tabular}{c|c|c|c}
 & 8-bus & 13-bus & 34-bus \\
\hline
Nominal & 1000 & 1000 & 1000 \\
\hline
Random & 968 & 925 & 932 
\end{tabular}
\caption{Number of times the relaxed problem is tight out of 1000 instances.}
\label{tab:sim}
\end{table}
As shown in Table \ref{tab:sim}, the relaxation is tight for all nominal situations. We offer some intuitive explanations for why this is the case. First consider the real power upper and lower bounds. Theorem \ref{thm:treelinebus} requires that when two buses are connected, not both have tight real power lower bounds. In the optimization problem we are minimizing the total system losses, so the feeder would try to meet the minimum power that is needed by the other nodes, since supplying more power will increase the total loss in the system. Therefore we expect that most of the buses to have $P_i^*= \ov{P}_i$. This is indeed the case in the simulations. Now consider the reactive power bounds. Theorem \ref{thm:treelinebus} requires that the lower reactive power bounds are not tight for the buses. In contrast to the real power, which flows downstream from the feeder to the end users, the reactive power flows up the tree from the end users to the feeder. This is because when the voltage is held constant, the users injected reactive power to support this voltage \cite{Oapos04}. Therefore for most of the nodes $Q_i^* >0$ in the simulation instances, so the lower bounds are not tight. In the random cases, since real power can flow up the tree, $Q_i$ could be positive or negative at bus $i$.      
\vspace{-0.4cm}
\section{Conclusion} \label{sec:V}
We studied the effects of constraints on power flow in a network and considered the
implication to the optimal power flow problem. We focused on the
injection region and showed how it can be used to understand the
optimal power flow problem. When there are no operation constraints,
we showed that the injection region is the entire upper half space.
 For tree networks, we showed that the injection region and its convex hull have the same Pareto-front when there is voltage magnitude constraints, line loss
constraints, line flow constraints, and some subset of bus power
constraints. 
%We also characterized the convex hull of some commonly
%occurring networks and a certain combination of them.

\section*{Acknowledgements}
This research was initiated while the second author was visiting the
Newton Institute in Cambridge, U.K., under the stochastic processes
in communication sciences program. Discussions with Frank Kelly in
the early stage of the research are much appreciated. Thanks also to
Alejandro Dominguez-Garcia and Javad Lavaei for comments on an
earlier version of this paper.

\begin{figure}[ht]
\centering
\includegraphics[scale=0.4]{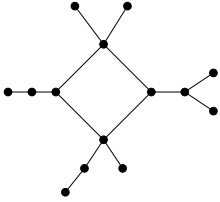}%
\caption{A distribution network with a ring feeder.}%
\label{fig:ring2}%
\end{figure}

\vspace{-0.4cm}
\appendix[Non-tree Networks] 
\label{sec:lossless}
Ideally, one would like to generalize the results for trees to
networks with cycles. However, this is difficult.
%%\footnote{The
%%result in \cite{Holst03} is actually stronger. It stated that if the
%%graph is not a tree, then there exists a positive semidefinite
%%matrix fitting the graph with rank $n-2$.}
We state some partial
results in this section, and they will be different than the result
stated in Theorem \ref{thm:treelinebus} in three aspects
\begin{itemize}
    \item We focus on lossless networks.
    \item Only voltage constraints are considered.
    \item We look at the convex hull instead of the Pareto-front.
\end{itemize}
Therefore the results in this section are of a weaker flavor than
Theorem \ref{thm:treelinebus} since we need to assume that the
networks are lossless and we only consider voltage constraints.
The results here are useful since in practice some distribution networks consists of a ring feeder and trees hanging of the feeder nodes as in Figure \ref{fig:ring2}. 
In this case, the objective functions are often to minimize the loss at the feeders. Also, the feeder nodes are generally considered as slack buses, so they only have a voltage constraint. 
Since minimizing a linear function over $\mathcal{A}$
and $\convhull(\mathcal{A})$ has the same objective values, characterizing the convex hull of the injection region is useful. 
%In practice, often the cost functions are taken to be linear in the
%power injections. If the convex hull of the injection region can be
%characterized, then the optimization with linear functions can be
%done efficiently.

The voltage constraint injection region is defined as
\begin{equation} \label{eqn:regionV}
\mathcal{P} = \{ \bd{p} : \bd{p}=\Real(\diag(\bd{v}\bd{v}^H \bd{Y}^H)), \ul{V}_i \leq |V|_i \leq \ov{V}_i \}.
\end{equation}
We can again define a enlarged convex region $\tl{\mathcal{P}}$ as
\begin{equation}
\tl{\mathcal{P}} = \{ \bd{p} : \bd{p}=\Real(\diag(\bd{W} \bd{Y}^H)), \ul{V}_i^2 \leq W_{ii} \leq \ov{V}_i^2,\bd{W} \sdp 0 \}. \label{eqn:regionW}
\end{equation}

We have the following theorem
\begin{thm} \label{thm:lossless}
Given a network with $n$ buses represented by its bus admittance matrix $\bd{Y}$. Let $\mathcal{P}$ and $\tl{\mathcal{P}}$ be defined as in \eqref{eqn:regionV} and \eqref{eqn:regionW} respectively. Then if the network is a {\em lossless cycle} or a {\em lossless cycle with one chord}, then
$\convhull(\mathcal{P})=\tl{\mathcal{P}}$.
\end{thm}

The next theorem states that joining the basic types of networks in a certain way preserves the characterization result. Given two networks $G$ and $H$, the network $K$ is said to be a 1-connection of $G$ and $H$ if it is possible to decompose $K$ into two components $K_1$ and $K_2$ such that they have only one node in common and no edges between them, where $K_1$ is equal to $G$ and $K_2$ is equal to $H$. Note by equal we mean that the admittance matrices are identical. In particular, if a line in $G$ or $H$ is lossless then its corresponding line in $K$ is also lossless. We say $K$ is obtained by 1-connecting $G$ and $H$. Figure \ref{fig:ring2} gives an example of a network obtained by 1-connecting a cycle and a number of trees.
\begin{thm} \label{thm:3}
Given a network on $n$ nodes with voltage constraints. Then $\convhull{\mathcal{P}} = \tl{\mathcal{P}}$ if the network is a result of repeatedly 1-connecting a lossless cycle and a tree.
\end{thm}

It is simple to check if a network has the topology that satisfies the conditions in Theorem \ref{thm:3}. Given a network, first decompose it into its one connected parts which can be done in linear time. Then one simply check each of the parts to see if they are a tree or a lossless cycle. 
%\begin{figure}[ht]
%\centering
%\subfigure[]{
%\includegraphics[scale=0.3]{5tree.eps}
%}
%\subfigure[]{
%\includegraphics[scale=0.3]{5cycle.eps}
%}
%\subfigure[]{
%\includegraphics[scale=0.3]{5cyclechord.eps}
%}
%\caption{Three networks: (a) a tree, (b) a cycle and (c) a cycle with a chord.}
%\label{fig:examples}
%\end{figure}
%\begin{figure}[ht]
%\centering
%\includegraphics[scale=0.2]{1connected.eps}
%\caption{A graph resulting from 1-connecting the basic components.}
%\label{fig:1connected}
%\end{figure}

%We suspect that a much stronger statement can be made about lossless
%networks. In particular, suppose $Y$ is a lossless network and is
%one of the two topologies in Theorem \ref{thm:lossless}, then we
%conjecture
%\begin{align*}
%& \{\bd{p} : \bd{p} =\Real{\diag(\bd{v}\bd{v}^H \bd{Y}^H)}, \ul{V}_i \leq |V|_i \leq \ov{V}_i \; \forall i\} \\
%= &\{\bd{p}: \bd{p} =\Real{\diag(\bd{v}\bd{v}^H \bd{Y}^H)}, |V|_i \leq \ov{V}_i \; \forall i\} \\
%= &\{\bd{p}: \bd{p} =\Real{\diag(\bd{v}\bd{v}^H \bd{Y}^H)}, |V|_i = \ov{V}_i \; \forall i\}.
%\end{align*}
%The intuition is that in a lossless network, the power flow is
%driven by the angle difference between the buses, so every injection
%vector should be achievable by adjusting the phases of the complex
%voltage vectors while keeping the magnitudes constant and there
%would be no holes in the injection region.

First we prove Theorem \ref{thm:lossless}. This requires that we prove an analogous result about trees first. Consider the following lemma
\begin{lem} \label{lem:hulltree}
Given a {\em tree} network with $n$ buses.  Let $\mathcal{P}$ and $\mathcal{P}_{\bd{W}}$ be defined as in \eqref{eqn:regionV} and \eqref{eqn:regionW} respectively. Then $\convhull (\mathcal{P})= \mathcal{P}_{W}$.
\end{lem}
\begin{IEEEproof}
To prove this theorem, it suffices to prove that minimizing linear functions over $\mathcal{P}$ and $\mathcal{P}_{\bd{W}}$ has the same optimal objective value for all coefficients \cite{Boyd04}. So consider the optimization problem
\begin{align}\label{eqn:opfPV}
J= \mbox{minimize } & \sum_{i=1}^n c_i P_i\\
 \mbox{subject to } & \ul{V}_i \leq |V_i| \leq \ov{V}_i, \; \forall i \nn \\
&  \bd{p}=\Real(\diag(\bd{v}\bd{v}^H \bd{Y}^H)), \nn
\end{align}
and its relaxation
\begin{align} \label{eqn:opfrW}
J_1= \mbox{minimize } & \sum_{i=1}^n c_i P_i \\
 \mbox{subject to } & \ul{V}_i^2 \leq W_{ii} \leq \ov{V}_i^2, \; \forall i \nn \\
 & \bd{W} \sdp 0 \nn \\
 & \bd{p}=\Real(\diag(\bd{W} \bd{Y}^H)), \nn
\end{align}
where the costs can be general (no longer constraint to be positive). We show $J=J_1$ for all $c_i$'s.

The dual of \eqref{eqn:opfrW} is
\begin{align} \label{eqn:opfD}
J_1=\mbox{maximize } & \sum_{i=1}^n (\ul{V}_i^2 \ul{\la}_i-\ov{V}_i^2 \ov{\la}_i) \\
 \mbox{subject to } & \bd{\La} + \bd{M} \sdp 0. \nn
\end{align}
From Lemma \ref{lem:treefit} if the costs are such that $\bd{M}$ is connected, then the optimal solution to \eqref{eqn:opfrW} is rank $1$ and clearly $J=J_1$. If $\bd{M}$ is disconnected, then $\bd{M}$ can be written as a block diagonal matrix. Suppose there are $K$ connected components of $\bd{M}$, then $\bd{M}=\diag(\bd{M}_1,\dots,\bd{M}_K)$. Since the network is a tree, $\bd{M}_i$ fits the topology of a tree for each $i$. Then \eqref{eqn:opfrW} and \eqref{eqn:opfD} decomposes into $K$ independent primal-dual subproblems, and we may apply Lemma \ref{lem:treefit} to each of them. Let $W_1^*, \dots,W_K^*$ denote the optimal solutions to each of the subproblems. By Lemma \ref{lem:treefit}, they are all rank $1$ so we can write $W_i^*=\bd{v}_i^* (\bd{v}_i^*)^H$ for each $i$. An optimal solution $\bd{W}^*$ to the original problem is given by $\bd{W}^*=\bd{v}^* (\bd{v}^*)^H$ where $\bd{v}^*=\bma \bd{v}_1^* \\ \vdots \\ \bd{v}_K^*\ebma$.
\end{IEEEproof}

Now we prove Theorem \ref{thm:lossless}. The approach is the same as in the proof of Lemma \ref{lem:hulltree}. That is, we look at \eqref{eqn:opfPV}, \eqref{eqn:opfrW} and their dual \eqref{eqn:opfD}. We say a matrix $\bd{A}$ is lossless if all the off diagonal terms of $A$ are purely imaginary or $0$. We prove the following lemma
\begin{lem} \label{lem:lossless}
Given a graph on $n$ nodes that is either an odd cycle or a cycle with one chord, if $\bd{A}$ is lossless, positive semidefinite and fits $G$, then $\rank(\bd{A}) \geq n-1$.
\end{lem}
Theorem \ref{thm:lossless} can be proved from Lemma \ref{lem:lossless}.  Suppose the electrical network is lossless and has the topology of an odd cycle or a cycle with one chord. The network being lossless means $\bd{Y}$ is purely imaginary, and $\bd{M}=\frac12 (\bd{C}\bd{Y}+\bd{Y}^H \bd{C})$ is also purely imaginary since $\bd{C}$ is real. Suppose that the costs are such that $M_{ik} \neq 0$ if $(i,k)$ is connected by a line in the network. Since $\bd{\La}^*$ is diagonal, the dual matrix $\bd{\La}^*+\bd{M}$ is positive semidefinite, lossless and fits the network topology. Apply Lemma \ref{lem:lossless} shows $\bd{W}^*$ is rank 1. If the cycle is even, we add a chord between two buses, and let the admittance of that chord go to $0$. Since all the functions in \eqref{eqn:opfrW} are continuous, the optimal solution of the network with a chord approaches the network without the chord as the admittance goes to $0$.

If the costs are such that $M_{ik}=0$ even if $(i,k)$ is connected in the network, then $\bd{M}$ either fits a tree or becomes disconnected. If $\bd{M}$ fits a tree, then apply Lemma \ref{lem:treefit}. If $\bd{M}$ becomes disconnected, then $\bd{M}$ can be written as a block diagonal matrix. If there are $K$ connected components of $\bd{M}$, then $\bd{M}=\diag(\bd{M}_1,\dots,\bd{M}_K)$. Since the network is a cycle (with a chord), then $\bd{M}_i$ is either a tree or a cycle for each $i$. We can apply Lemma \ref{lem:treefit} or Lemma \ref{lem:lossless} to each component and obtain an optimal solution $\bd{W}^*$ in the same way as in the tree network case. To finish Theorem \ref{thm:lossless}, it remains to proof Lemma \ref{lem:lossless}.
\begin{IEEEproof}
Given a graph $G$, the tree-width of $G$ is a number that intuitively captures how close $G$ is to a tree. For example, the tree-width of a tree is $1$, and the tree-width of a cycle is $2$. The rigorous definition and some methods of computing the tree-width  the reader may consult \cite{Bodlaender96}. A graph of tree-width 2 is also called serial-parallel graph or a partial-2-tree. The following lemma collects the known results that we need.
\begin{lem}\label{lem:fit}
If $G$ is a cycle of length $n$, then the minimum rank of real positive semidefinite matrices fitting $G$ is $n-2$  \cite{Fallat07}. More generally, if the graph has tree-width $2$, the minimum rank is $n-2$ \cite{Holst03, Johnson09}.
\end{lem}
Given a graph $G$ with $n$ nodes and $m$ edges. We construct a bipartite graph derived from $G$ that we call the bipartite expansion of $G$ and denote by $B(G)$.  $B(G)$ is a bipartite graph with $2n$ nodes and $2m$ edges. Label the nodes $1,2,\dots,n,1',2',\dots,n'$ with the bipartition being $\{1,\dots,n\}$ and $\{1',\dots,n'\}$. There is an edge between $i$ and $k'$ if and only if $i \neq k$ and $(i,k)$ is an edge in $G$. If $G$ is an odd cycle then $B(G)$ is also a cycle and if $G$ is a cycle with a chord then $B(G)$ has tree-width 2 (a subclass of linear-2-trees in the language of \cite{Johnson09}). Two examples are given in Figures \ref{fig:3BG} and \ref{fig:4BG}. If $G$ is an even cycle then $B(G)$ is two disconnected cycles, therefore the assumption of odd cycle is needed in the Lemma.
\begin{figure}[ht]
\centering
\subfigure[$G$]{
    \includegraphics[scale=0.2]{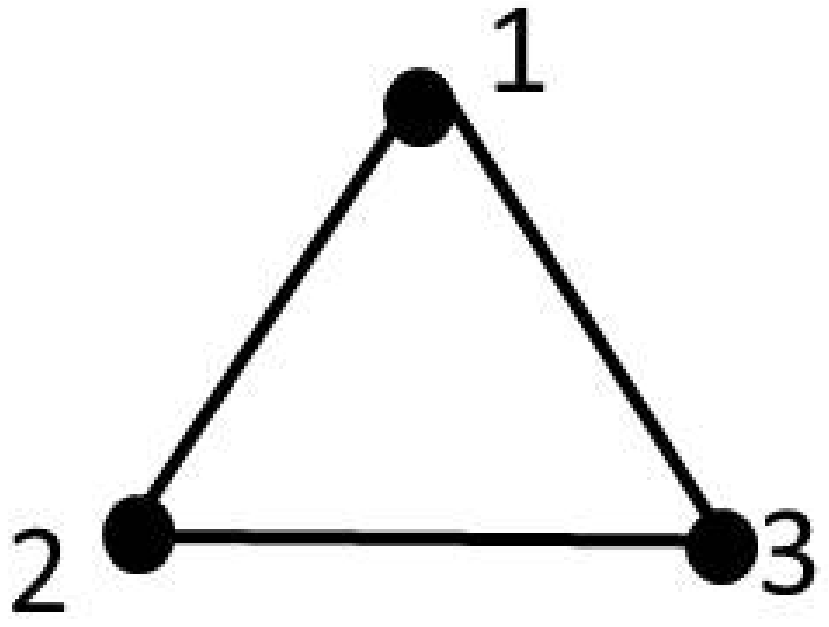}
}
\subfigure[$B(G)$]{
    \includegraphics[scale=0.2]{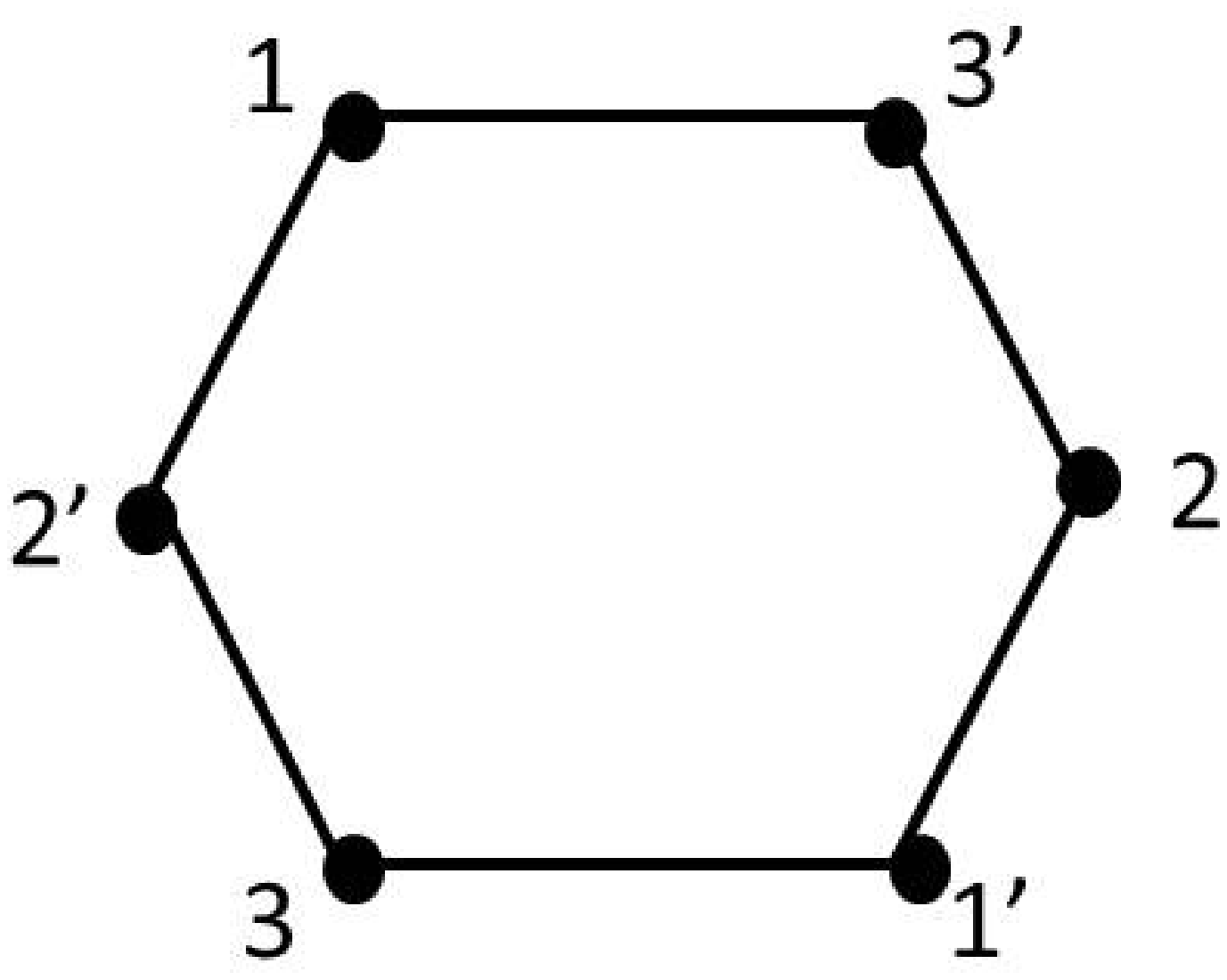}
}
\caption{(a) shows a 3-cycle and (b) shows its bipartite expansion.}
\label{fig:3BG}
\end{figure}
\begin{figure}[ht]
\centering
\subfigure[$G$]{
    \includegraphics[scale=0.2]{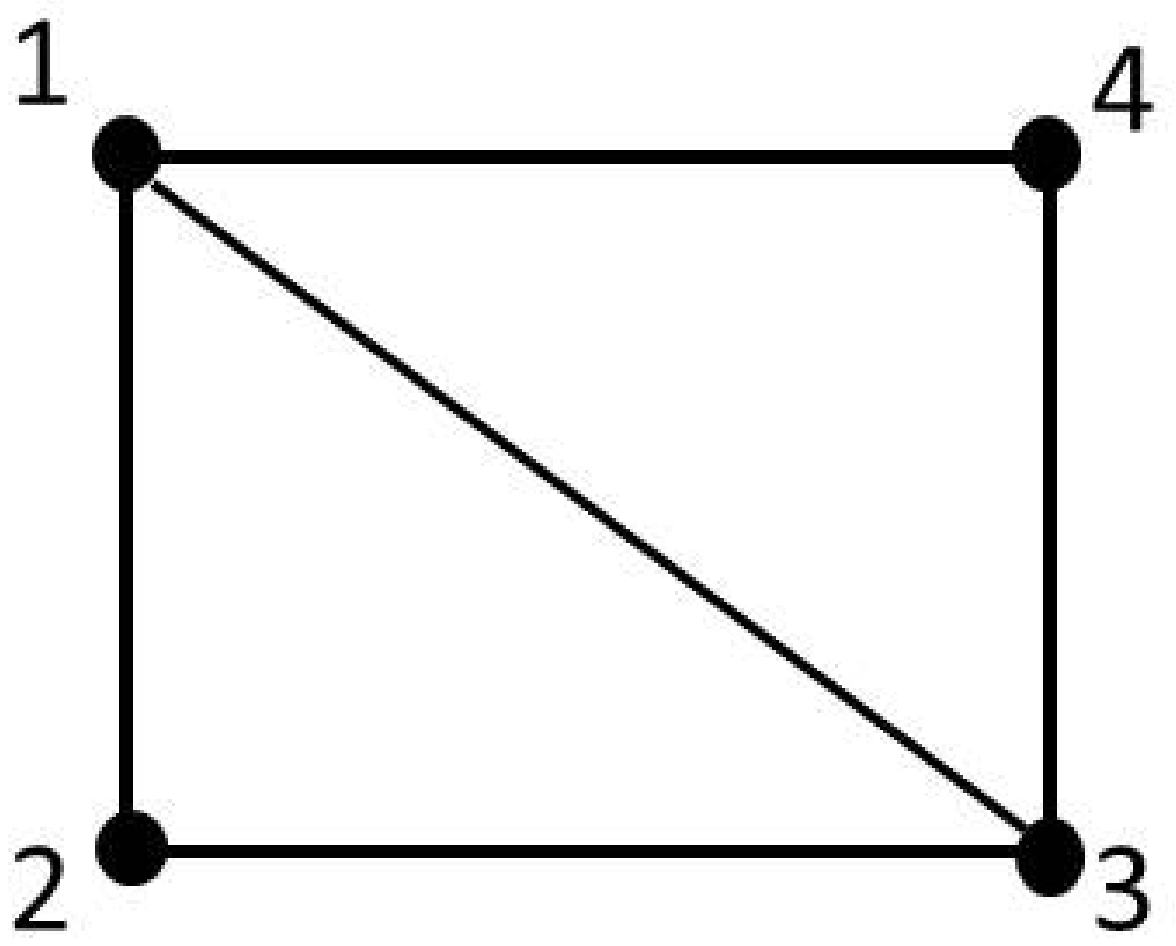}
}
\subfigure[$B(G)$]{
    \includegraphics[scale=0.2]{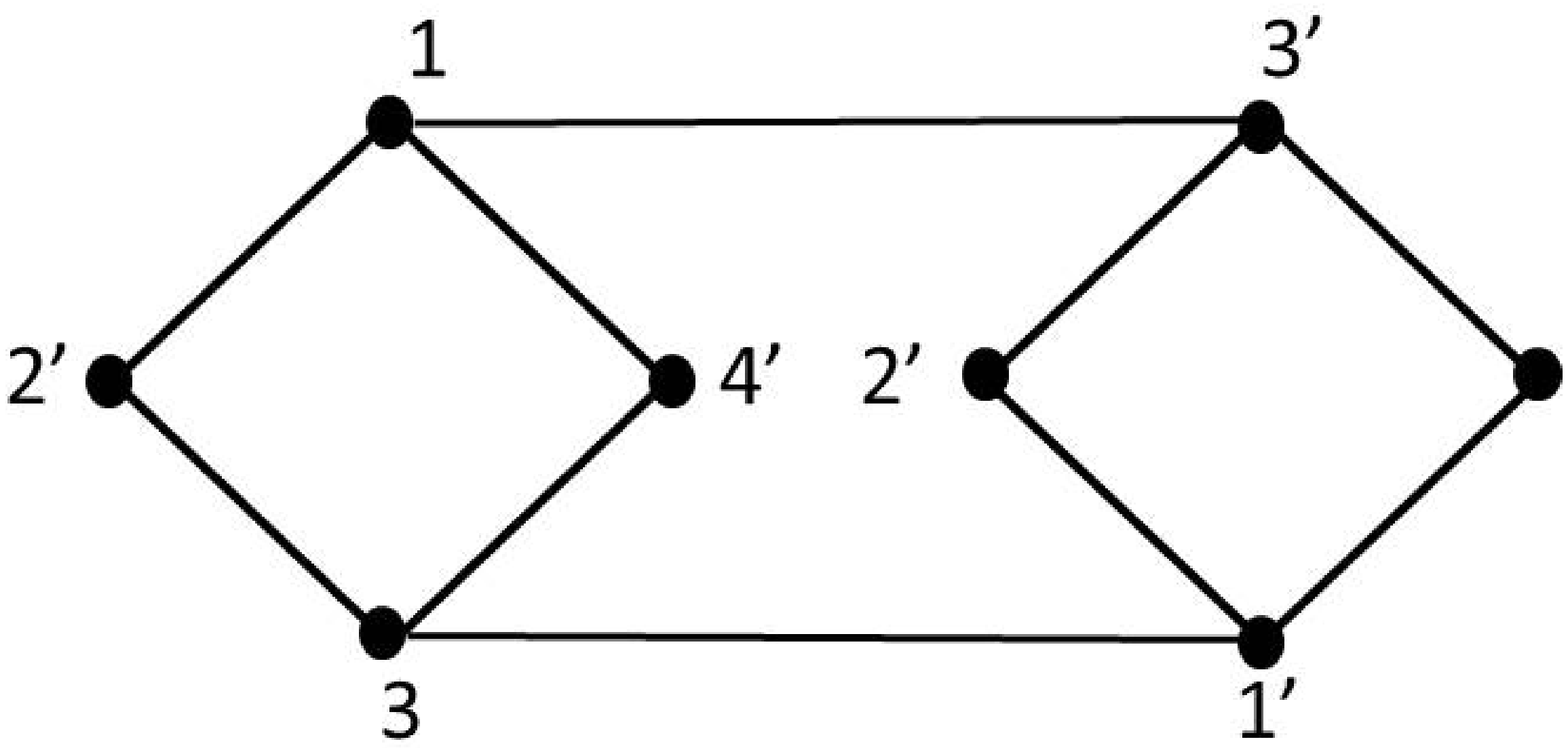}
}
\caption{(a) shows a 4-cycle with a chord and (b) shows its bipartite expansion.}
\label{fig:4BG}
\end{figure}

Given a graph $G$, suppose $\bd{A}$ is lossless, positive semidefinite and $\bd{A}$ fits $G$. We show that the rank of $\bd{A}$ cannot be lower than $n-1$. Suppose $\bd{A}$ has rank $r$. Then $\bd{A}$ can be factored as $\bd{A}=\bd{Z}^H \bd{Z}$ for some complex matrix $r \times n$ matrix $\bd{Z}$. Let $\z_1,\dots,\z_n \in \C^r$ be the columns of $\bd{Z}$. They satisfy the graph topology condition
\begin{equation}\label{eqn:topo}
\z_i^H \z_k =
\begin{cases}
0 &\mbox{ if } i \nsim k \\
\neq 0 &\mbox{ if } i \sim k
\end{cases}
\end{equation}
and the lossless line condition
\begin{equation} \label{eqn:lossless}
\Real(\z_i^H \z_k) =0 \mbox{ if } i \neq k.
\end{equation}
From each complex vector we define two real vectors as
\begin{equation*}
\x_i = \bma \Real(\z_i) \\ \Imag(\z_i) \ebma \;\;\;\;\; \y_i = \bma \Imag(\z_i) \\ -\Real(\z_i) \ebma
\end{equation*}
Since $\z_i \in \C^r$, then $\x_i, \y_i \in \R^{2r}$. By algebra, $\Real(\z_i^H \z_k)=\x_i^T \x_k=\y_i^t \y_k$ and $\Imag(\z_i^H \z_k)=\x_i^T \y_k$. In terms of $\x$'s and $\y$'s, \eqref{eqn:topo} becomes
\begin{equation}\label{eqn:topox}
\x_i^T \y_k =
\begin{cases}
0 &\mbox{ if } i \nsim k \\
\neq 0 &\mbox{ if } i \sim k
\end{cases}
\end{equation}
and \eqref{eqn:lossless} becomes
\begin{equation}\label{eqn:losslessx}
\x_i^T \x_k = \y_i^T \y_k =0 \mbox{ if } i \neq k.
\end{equation}
Define the matrix $B$ to be the $2r \times 2n$ matrix with columns $\x_1,\dots,\x_n,\y_1,\dots,\y_n$. By \eqref{eqn:topox} and \eqref{eqn:losslessx} $B$ fits $B(G)$. But if $G$ is an odd cycle or a cycle with one chord, applying Lemma \ref{lem:fit} to $B(G)$ gives $\rank(B) \geq 2n-2$. Thus $2r \geq 2n-2$ or $r \geq n-1$.
\end{IEEEproof}

Now we proceed to the proof Theorem \ref{thm:3}. Given a network $G$, we say the matrix $\bd{A}$ satisfies $G$ if $\bd{A}$ fits the topology of $G$ and $A_{ik}$ is purely imaginary if the line from bus $i$ to bus $k$ is lossless. We have the following lemma.
\begin{lem}\label{lem:1connected}
Given two networks $G$ and $H$ with $n$ and $m$ buses respectively, let $K$ be a network obtained by 1-connecting $G$ and $H$, so $K$ has $n+m-1$ buses. If $\bd{A}$ is a positive semidefinite matrix that satisfies $K$, then $\rank(\bd{A}) \geq n+m-2$.
\end{lem}
From the basic topologies in Theorem \ref{thm:lossless}, we can apply the Lemma \ref{lem:1connected} repeatedly to get Theorem \ref{thm:3}. A version of Lemma \ref{lem:1connected} just about graphs (without considering lossless lines and such) is known in the graph theory community \cite{Fallat07,Beagley07}. We give a proof here to show the additional condition of lossless lines does not change the result.
\begin{IEEEproof}
Let $G$, $H$ and $K$ be networks given in the statement of the Lemma.  Label the buses in $K$ to be $1,2,\dots,n-1,n,n+1,n+2,\dots,n+m-1$ where the subnetwork induced by $1,\dots,n-1,n$ corresponds to $G$ and the subnetwork induced by $n,n+1,n+m-1$ corresponds to $H$. So bus $n$ is the common bus in the 1-connection. Suppose $\bd{A}$ is a $(n+m-1) \times (n+m-1)$ positive semidefinite matrix that satisfies $K$ and has rank $r$. Then it is possible to factor $\bd{A}$ as $\bd{A}=\bd{Z}^H \bd{Z}$ for some $r \times (n+m-1)$ matrix $\bd{Z}$. Let $\z_1,\dots,\z_{n+m-1}$ be the columns of $\bd{Z}$. Let $\mathcal{U}$ be the subspace spanned by $\z_1,\dots,\z_{n-1}$ and $\mathcal{V}$ be the subspace spanned by $\z_{n+1}, \dots, \z_{n+m-1}$. By construction of $K$, there are no lines between the set of buses $\{1,\dots,n-1\}$ and $\{n+1,\dots,n+m-1\}$. Therefore $\mathcal{V}$ is orthogonal to $ \mathcal{U}$. We may write vector $\z_n$ as $\z_n=\bd{u}+\bd{v}+\bd{w}$ where $\bd{u} \in \mathcal{U}$, $\bd{v} \in \mathcal{V}$ and $\bd{w}$ is orthogonal to $\mathcal{U}$ and $\mathcal{V}$. Let $\bd{Z}_G$ be the matrix with columns $\z_1,\dots,\z_{n-1},\bd{u}$ and $\bd{Z}_H$ be the matrix with columns $\bd{v},\z_{n+1},\dots,\z_{n+m-1}$. Let $\bd{A}_G=\bd{Z}_G^H \bd{Z}_G$. Since $\z_i^H \bd{u} = \z_i^H \z_n$ for $i=1,\dots,n-1$, $\bd{A}_G$ equals the matrix formed by the first $n$ rows and $n$ columns of $\bd{A}$. By the assumption $\bd{A}$ satisfies $K$, so $\bd{A}_G$ satisfies $G$. Similarly $\bd{Z}_H^H \bd{Z}_H$ satisfies $H$. By the assumption in the Lemma, we have $\rank(\bd{Z}_G) \geq n-1$ and $\rank(\bd{Z}_H) \geq m-1$, so equivalently $\dim \mathcal{U} \geq n-1$ and $\dim \mathcal{V} \geq m-1$. Since $\mathcal{U}$ is orthogonal to $\mathcal{V}$ and $\z_1,\dots,\z_{n+m-1}$ spans $\mathcal{U} + \mathcal{V}$, $\rank(\bd{A}) = \dim \mathcal{U}+\dim \mathcal{V} \geq (n-1)+(m-1)=n+m-2$.
\end{IEEEproof}

\appendix[Proof of Theorem \ref{thm:1}]
The following basic lemma from linear algebra is useful.
\begin{lem}[Rank Nullity Theorem]\label{lem:rangeA}
Let $A$ be a $n \times n$ real symmetric matrix. Let $\im(A)$ and $\ker(A)$ denote the image and kernel of $A$, respectively. Then
$\dim \im(A) +\dim \ker(A) = n$ and
$\im(A) \oplus \ker(A)=\R^n$,
where $\oplus$ is the direct sum.
\end{lem}

First consider the case where the network is lossless. Then any feasible injection vector must be on the conservation of energy plane. We need to show that any point on the plane can be achieved. Since the network is lossless $\bd{Y}=j \Imag(\bd{Y})$ where $\Imag(\bd{Y})$ is a $n \times n$ real symmetric matrix and each row of $\Imag(\bd{Y})$ sums to $0$ by \eqref{eqn:Y}. Therefore $\Imag(\bd{Y})$ is a generalized graph Laplacian matrix where the admittances can be interpreted as weights on the edges. By a standard result in graph theory, $\dim \ker(\Imag(\bd{Y})) =1$ and $\ker (\Imag(\bd{Y}))$ is spanned by the all one's vector $\bd{1}$. By Lemma \ref{lem:rangeA}, $\im (\Imag(\bd{Y}))$ is the linear subspace in $\R^n$ orthogonal to $\bd{1}$. Let $\bd{p}^0$ be an injection vector on the conservation of energy plane, that is $\sum_{i=1}^n P_i^0 =0$. Since $\bd{1}^T \bd{p}^0 =0$, there is a unique vector $\bd{v}^0$ such that $\bd{Y} \bd{v}^0 =\bd{p}^0$ and $\bd{1}^T \bd{v}^0 =0$. Choose the voltage vector $\bd{v}=(-\bd{v}^0+j \bd{1})$, then
\begin{align}
& \Real(\diag((-\bd{v}^0+j\bd{1})(-\bd{v}^0+j\bd{1})^H \bd{Y}^H)) \\
&= \Real(\diag((\bd{v}^0 \bd{1}^T + \bd{1} (\bd{v}^0)^T)\Imag(\bd{Y})) \\
&\mbox{ }+j \diag((\bd{v}^0 (\bd{v}^0)^T+\bd{1}\bd{1}^T)\Imag(\bd{Y}))) \nn \\
&\stackrel{(a)}{=} \bd{p}^0,\nn
\end{align}
where $(a)$ follows from the choice of $\bd{v}^0$ and $\Imag(\bd{Y})$ being symmetric. This finishes the proof for a lossless network.

Next consider the case where the network is lossy. The proof proceeds in two parts, first we show that the conservation of energy boundary $\sum_{i=1}^n P_i=0$ can be arbitrarily closely from above, and then we show the injection region is convex. Since the network is lossy, $\Real(\bd{Y})$ is a $n \times n$ real positive semidefinite Laplacian matrix. By conservation of energy, any power injection vector achieved must satisfy $\sum_{i=1}^n P_i >0$ if $\bd{p} \neq 0$. Let $\bd{p}^0$ be a vector on the conservation of energy plane. We show there is a voltage vector $\bd{v}$ that achieves a point arbitrarily close to $\bd{p}^0$. Since $\bd{1}^T \bd{p}^0 =0$, by Lemma \ref{lem:rangeA} there is a unique vector $\bd{v}^0$ such that $\Real(\bd{Y}) \bd{v}^0 =\bd{p}^0$ and $\bd{1}^T \bd{v}^0 =0$. Let $\bd{v}=(\al \bd{1} +\frac{1}{\al} \bd{v}^0)$ for some $\al \geq 0$ and the corresponding injection vector $\bd{p}$ is
\begin{align}
\bd{p} &=\Real(\diag(\bd{v}\bd{v}^T \bd{Y})) \\
&= \Real(\diag ((\al \bd{1} + \frac1\al \bd{v}^0)(\al \bd{1} + \frac1\al \bd{v}^0)^T (\Real(\bd{Y})+j \Imag(\bd{Y}))) \nn \\
&= \diag ((\al^1 \bd{1}\bd{1}^T + \bd{v}^0 \bd{1}^T+ \bd{1} (\bd{v}^0)^T +\frac{1}{\al^2}\bd{v}^0 (\bd{v}^0)^T) \Real(\bd{Y})) \nn \\
& \stackrel{(a)}{=} \diag(\bd{1} (\bd{v}^0)^T \Real(\bd{Y})) + \frac{1}{\al^2} \diag(\bd{v}^0 (\bd{v}^0)^T \Real(\bd{Y})) \nn \\
&\stackrel{(b)}{=} \diag(\bd{1} (\bd{p}^0)^T)+\frac{1}{\al^2} \diag(\bd{v}^0 (\bd{p}^0)^T) \nn \\
&= \bd{p}^0 + \frac{1}{\al^2} \diag(\bd{v}^0 (\bd{p}^0)^T) \nn,
\end{align}
where $(a)$ follows from $\bd{1} \in \ker (\Real(\bd{Y}))$ and $\Real(\bd{Y})$ is symmetrical, $(b)$ follows from the choice of $\bd{v}^0$. We can increase $\al$ to make $\bd{p}$ arbitrarily close to $\bd{p}^0$. For example, if we want $||\bd{p}-\bd{p}^0||_\infty \leq \eps$, then choose
\begin{equation*}
\al \geq \sqrt{\frac{||\bd{p}^0||_\infty ||\bd{v}^0||_\infty}{\eps}}.
\end{equation*}

The next lemma states that $\mathcal{P}$ is convex.
\begin{lem}\label{lem:cone}
The injection region $\mathcal{P}$ as defined in eqn. (\ref{eq:energy}) is a convex set.
\end{lem}
Theorem \ref{thm:1} follows from Lemma \ref{lem:cone}. Since the injection region is convex, and the boundary $\sum_{i=1}^n P_i=0$ can be approached arbitrarily closely from above, it includes the open half upper space. In addition the origin can be achieved using the all zeros voltage vector. It remains to prove the lemma.
\begin{IEEEproof}
For a given network with $n$ buses represented by $\bd{Y}$, define $\mathcal{P}_{\ov{V}}$ as
\begin{equation}\label{eqn:regionV2}
\mathcal{P}_{\ov{V}}=\{\bd{p}\in \R^n : \bd{p}=\Real(\diag(\bd{v}\bd{v}^H\bd{Y}^H)), ||\bd{v}||_2 \leq \ov{V}\},
\end{equation}
where $||\bd{v}||_2=(\sum_{i=1}^n |V|_i^2)^{\frac12}$. $\mathcal{P}_{\ov{V}}$ approaches the unconstrained injection region as $\ov{V}$ tends to infinity.  $\mathcal{P}_{\ov{V}}$ cannot have holes since if $\bd{p} \in \mathcal{P}_{\ov{V}}$, then $\al \bd{p} \in \mathcal{P}_{\ov{V}}$ for $\al \in [0,1]$. Therefore to prove the convexity of $\mathcal{P}_{\ov{V}}$ it suffices to prove it has convex boundary. Consider the optimization problem
\begin{align}\label{eqn:opfS}
J=\mbox{minimize } &\sum_{i=1}^n c_i P_i \\
\mbox{subject to } &||\bd{v}||_2 \leq \ov{V} \nn \\
& \bd{p}=\Real(\diag(\bd{v}\bd{v}^HY^H)). \nn
\end{align}
Relaxing and eliminating $\bd{p}$, we get
\begin{align}\label{eqn:opfSW}
J_1=\mbox{minimize } &\Tr(\bd{M} \bd{W}) \\
\mbox{subject to } & \sum_{i=1} W_{ii} \leq \ov{V}^2 \nn \\
& \bd{W} \sdp 0, \nn
\end{align}
By changing the costs, we are exploring the boundaries of the two regions with linear functions. We want to show that all the point on the boundary of the larger region is in fact in the smaller region.

First we show that for all $\bd{M}$ there is an optimal $\bd{W}^*$ for \eqref{eqn:opfSW} which is rank 1.  To solve \eqref{eqn:opfSW}, expand $\bd{W}$ in terms of its eigenvectors, so $\bd{W}=w_1 \w_1 \w_1^H + \cdots w_n \w_n \w_n^H$ where $\w_i$ is unit norm and $\sum_{i=1}^n w_i \leq \ov{V}^2$. Then \eqref{eqn:opfSW} can be written as
\begin{align}\label{eqn:opfmW}
\mbox{minimize } & \sum_{i=1}^n w_i \w_i^H M\w_i \\
\mbox{subject to } & \sum_{i=1}^n w_i \leq \ov{V}^2 \nn \\
& \bd{W}=\sum_{i=1}^n (w_i \w_i \w_i^H) \sdp 0. \nn
\end{align}
By the well known result about Rayleigh quotients \cite{Horn85}, to minimize any of the terms $\w_i^H M \w_i$, the optimal $\w_i^*=\bd{m}_1$, where $\bd{m}_1$ is the eigenvector corresponding to the smallest eigenvector of $M$. Therefore the optimal solution to \eqref{eqn:opfSW} is $\bd{W}=\sum_{i=1}^n w_i \bd{m}_1 \bd{m}_1^H= \ov{V}^2 \bd{m}_1 \bd{m}_1^H$ and is rank 1.

If $\bd{m}_1$ is not unique, since eigenvector are not continuous in the entries of the matrix, we can perturb $\bd{Y}$ by an arbitrarily small amount to obtain a $M$ that has a unique eigenvector corresponding to the smallest value. Note the power vector $\bd{p}$ is continuous in the entries of $Y$. From uniqueness of $\bd{m}_1$ and the fact there is no gap between \eqref{eqn:opfS} and \eqref{eqn:opfSW}, the two regions have the same boundary. Taking $\ov{V}$ to infinity finishes the proof.
\end{IEEEproof}

\bibliography{mybib}
\bibliographystyle{IEEEtran}
\end{document}

%% file: header.tex
\usepackage{cite}
\usepackage{amsmath, amsthm, amssymb,graphicx,bbm}
\usepackage{verbatim}
\usepackage{psfrag}
\usepackage{algorithm2e}
\usepackage{bm}
\usepackage{color}
\usepackage[tight,footnotesize]{subfigure}

\newcommand{\R}{\mathbb{R}}

\newcommand{\al}{\alpha}

\newcommand{\bd}{\mathbf}

\newcommand{\tl}{\tilde}

\newcommand{\w}{\bd{w}}
\newcommand{\x}{\bd{x}}
\newcommand{\y}{\bd{y}}

\newcommand{\ov}{\overline}
\newcommand{\ul}{\underline}
\newcommand{\la}{\lambda}
\newcommand{\z}{\bd{z}}
\newcommand{\C}{\mathbb{C}}

\newcommand{\sdp}{\succcurlyeq}
\newcommand{\bma}{\begin{bmatrix}}
\newcommand{\ebma}{\end{bmatrix}}

\newcommand{\eps}{\epsilon}

\newcommand{\nn}{\nonumber}

\newcommand{\La}{\Lambda}

\newcommand{\mc}{\mathcal}

\DeclareMathOperator{\rank}{rank}
\DeclareMathOperator{\Tr}{Tr}

\DeclareMathOperator{\diag}{diag}
\DeclareMathOperator{\Real}{Re}
\DeclareMathOperator{\Imag}{Im}

\DeclareMathOperator{\im}{image}

\DeclareMathOperator{\convhull}{convhull}
\DeclareMathOperator{\conj}{conj}
\newtheorem{defn}{Definition}
\newtheorem{thm}{Theorem}
\newtheorem{claim}[thm]{Claim}
\newtheorem{lem}[thm]{Lemma}